\documentclass[12pt,oneside,reqno]{amsart}
\usepackage{}
\usepackage{amssymb}
\usepackage[colorlinks=true, linkcolor=blue, citecolor=red]{hyperref}
\usepackage{bbm}
\usepackage{cite}
\usepackage{amsmath}
\usepackage{graphicx}
\usepackage{mathrsfs}
\usepackage{stmaryrd}
\usepackage{color}
\usepackage{soul}
\usepackage[dvipsnames]{xcolor}
\usepackage{amsfonts}
\usepackage{enumerate,amsmath,amssymb,amsthm}

\numberwithin{equation}{section}

\newcommand{\be}{\begin{eqnarray}}
\newcommand{\mE}{\end{eqnarray}}
\newcommand{\ce}{\begin{eqnarray*}}
\newcommand{\de}{\end{eqnarray*}}
\newtheorem{theorem}{Theorem}[section]
\newtheorem{lemma}[theorem]{Lemma}
\newtheorem{remark}[theorem]{Remark}
\newtheorem{definition}[theorem]{Definition}
\newtheorem{proposition}[theorem]{Proposition}
\newtheorem{example}[theorem]{Example}
\newtheorem{corollary}[theorem]{Corollary}

\def\eps{\varepsilon}

\def\p{\partial}

\def\[{{\Big[}}
\def\]{{\Big]}}
\def\<{{\langle}}
\def\>{{\rangle}}
\def\({{\Big(}}
\def\){{\Big)}}

\def\bx{{\mathbf{x}}}

\def\dif{{\mathord{{\rm d}}}}

\def\no{\nonumber}
\def\={&\!\!=\!\!&}
\def\bt{\begin{theorem}}
\def\et{\end{theorem}}
\def\bl{\begin{lemma}}
\def\el{\end{lemma}}
\def\br{\begin{remark}}
\def\er{\end{remark}}

\def\bd{\begin{definition}}
\def\ed{\end{definition}}
\def\bp{\begin{proposition}}
\def\ep{\end{proposition}}
\def\bc{\begin{corollary}}
\def\ec{\end{corollary}}
\def\bx{\begin{example}}
\def\ex{\end{example}}

\def\cD{{\mathcal D}}

\def\cL{{\mathcal L}}

\def\mC{{\mathbb C}}

\def\mE{{\mathbb E}}

\def\mN{{\mathbb N}}

\def\mP{{\mathbb P}}

\def\mR{{\mathbb R}}

\def\sI{{\mathscr I}}
\def\sJ{{\mathscr J}}

\def\sL{{\mathscr L}}

\def\sN{{\mathscr N}}

\def\sR{{\mathscr R}}

\def\sV{{\mathscr V}}

\def\sX{{\mathscr X}}
\def\sY{{\mathscr Y}}

\def\geq{\geqslant}
\def\leq{\leqslant}

\allowdisplaybreaks

\begin{document}

\title{Diffusion approximation for multi-scale stochastic reaction-diffusion equations}

\date{}

\author{ Longjie Xie\,\, and\,\, Li Yang}


\address{Longjie Xie:
	School of Mathematics and Statistics $\&$ Research Institute of Mathematical Science, Jiangsu Normal University,
	Xuzhou, Jiangsu 221000, P.R.China\\
	Email: longjiexie@jsnu.edu.cn
}

\address{Li Yang:
	School of Mathematics, Shandong University,
	Jinan, Shandong 250100, P.R.China\\
	Email: llyang@mail.sdu.edu.cn
}

\thanks{This work is supported  by
the NSFC (No. 12090011, 12071186, 11771187, 11931004).}

\begin{abstract}

In this paper, we study the diffusion approximation for singularly perturbed stochastic reaction-diffusion equation with a fast oscillating term. The asymptotic limit  for the original system is obtained, where an extra Gaussian term appears. Such a term is explicitly given in terms of the solution of  Poisson equation in Hilbert space. Moreover, we also obtain the  rate of convergence, and the convergence rate is shown  not to depend  on the regularity of the coefficients of the original  system with respect to the fast variable, which coincides with the intuition since the fast component has been totally homogenized out in the limit equation.

	\bigskip

	
	\vspace{2mm}
\end{abstract}

\maketitle


\section{Introduction}

Let $T>0$ and  $D=(0,L)\subset\mR$ be a bounded inverval. Consider the following fully coupled slow-fast stochastic reaction-diffusion equation with Dirichlet
boundary condition:
\begin{equation} \label{spde0}
\left\{ \begin{aligned}
&\dif X^{\eps}_t(\xi) =\Delta X^{\eps}_t(\xi)\dif t+f(X^{\eps}_t(\xi), Y^{\eps}_t(\xi))\dif t+\dif W_t^1(\xi),\\
&\dif Y^{\eps}_t(\xi) =\eps^{-1}\Delta Y^{\eps}_t(\xi)\dif t+\eps^{-1}g(X^{\eps}_t(\xi), Y^{\eps}_t(\xi))\dif t+\eps^{-1/2} \dif W_t^2(\xi),\\
&X^{\eps}_t(0)=X^{\eps}_t(L)=Y^{\eps}_t(0)=Y^{\eps}_t(L)=0, \qquad t\in(0,T],\\
&X^{\eps}_0(\xi)=x(\xi),\,\,Y^{\eps}_0(\xi)=y(\xi),\qquad\xi\in D,
\end{aligned} \right.
\end{equation}
where  $f,g: \mR^2\to \mR$ are measurable functions, $W^1_t$ and $W^2_t$ are mutually independent  $L^2(D)$-valued $Q_1$- and $Q_2$-Wiener processes, and the small parameter $0<\eps\ll 1$ represents the separation  of  time scales.
Such multi-scale system appears frequently in many real-world  dynamical systems such as combustion, epidemic propagation  and dynamics of populations (see \cite{Ku,PS}). In such a system,   $X_t^\eps$ is called the slow process which can be thought
of as the mathematical model for a phenomenon appearing at the natural time scale, while  $Y_t^\eps$ (with time order $1/\eps$) is referred as the fast motion  which can be interpreted as the fast environment.

\vspace{1mm}
Usually, system of the form (\ref{spde0}) is difficult to deal with due to the two widely
separated time scales and the cross interactions between the fast and slow modes. Thus a simplified
equation which governs the evolution of the system over a long time scale is highly
desirable and is quite important for applications.
To give precise result, it is convenient to look at the equation in the abstract Hilbert space $H:=L^2(D)$, where
the system (1.1) can be rewritten as the stochastic partial differential equation (SPDE for short)
\begin{equation} \label{spde00}
\left\{ \begin{aligned}
&\dif X^{\eps}_t \!=\!AX^{\eps}_t\dif t+F(X^{\eps}_t, Y^{\eps}_t)\dif t+\dif W_t^1,\qquad\qquad\qquad\quad X^{\eps}_0=x\in H,\\
&\dif Y^{\eps}_t =\eps^{-1}AY^{\eps}_t\dif t+\eps^{-1}G(X^{\eps}_t, Y^{\eps}_t)\dif t+\eps^{-1/2} \dif W_t^2,\qquad Y_0^\eps=y\in H,
\end{aligned} \right.
\end{equation}
with $A:\cD(A)\subset H\rightarrow H$ being an unbounded linear operator,
and $F, G$ are Nemytskii operators defined by
\begin{align}\label{FG}
F(x,y)(\xi):=f(x(\xi),y(\xi))\quad\text{and}\quad G(x,y)(\xi):=g(x(\xi),y(\xi)).
\end{align}
Then the celebrated theory of averaging principle says that a good approximation
of the slow component in system (\ref{spde00}) can be obtained by averaging the coefficient with respect to
parameters in the fast variable. More precisely, under certain regularity assumptions
on the coefficients, the slow process  $X_t^\eps$ will converge  as $\eps\to 0$ to the solution of the following so-called averaged equation:
\begin{align} \label{lime}
\dif \bar X_t  =A\bar X_t\dif t+\bar F(\bar X_t)\dif t+\dif W_t^1,
\end{align}
where
$$
\bar F(x):=\int_{H}F(x,y)\mu^{x}(dy),
 $$
and $\mu^{x}(\dif y)$ is the unique invariant measure of the transition semigroup for the  frozen equation
\begin{align}\label{froz}
\dif Y_t^x=AY_t^x \dif t+G(x,Y_t^x)\dif t+\dif W_t^2, \quad Y_0^x=y\in H.
\end{align}
The reduced system (\ref{lime}) then captures the  essential dynamics of the original  system (\ref{spde00}), which does not depend on the fast variable any more and thus is much simpler.

\vspace{1mm}

In the past decades, the averaging principle for slow-fast systems has been intensively studied. We refer the readers to the fundamental paper  by Khasminskii \cite{K1} for stochastic differential equations (SDEs for short), see also \cite{BK,GR,HLi,KY,PXW,V0}. Generalization to the infinite dimensional setting is more difficult and  has been carried out only relative  recently. In \cite{CF}, Cerrai and Freidlin studied the averaging principle for a class of stochastic reaction-diffusion equations whose additive noise is included only in the fast motion. Later, Cerrai \cite{Ce,C2} extended the result in \cite{CF} to more general cases, see also \cite{Br1,DS,FWL0,GP1,GP2,GP3,LRS1,LW,PXG,WR} and the references therein for further developments.

 \vspace{1mm}
In this paper, we consider the following fully coupled multi-scale stochastic reaction-diffusion equation in the Hilbert space $H$:
\begin{equation} \label{spde11}
\left\{ \begin{aligned}
&\dif X^{\eps}_t \!=\!AX^{\eps}_t\dif t+F(X^{\eps}_t, Y^{\eps}_t)\dif t+\eps^{-1/2}B(X^{\eps}_t, Y^{\eps}_t)\dif t+\Sigma(X^{\eps}_t, Y^{\eps}_t)\dif W_t^1,\,\,\,\, X^{\eps}_0=x,\\
&\dif Y^{\eps}_t =\eps^{-1}AY^{\eps}_t\dif t+\eps^{-1}G(X^{\eps}_t, Y^{\eps}_t)\dif t+\eps^{-1/2} \dif W_t^2,\qquad\qquad\qquad\qquad\quad\! Y^{\eps}_0=y.
\end{aligned} \right.
\end{equation}
Compared with the system (\ref{spde00}) and all the above-mentioned papers, the main feature of  SPDE (\ref{spde11}) is that even the slow process $X_t^\eps$ has a fast varying component.  This is known to be  important   for applications in homogenization, which has its own interest in the theory of PDEs, see e.g. \cite{HP, HP2} and \cite[Chapter IV]{Fr}.  Moreover, such singularly perturbed equation (with the appearance of a fast term in the slow equation) provides a framework to model many physical systems, from colloidal particles in a fluid \cite{MJ,Nel} to a camera tracking an object \cite{Pap}. We refer the interested readers to   \cite[Section 11.7]{PS} for more applications. In fact, a very particular case of the equation (\ref{spde11}) is the following   Langevin equation:
\begin{align}\label{lav}
\eps\ddot{X_t^\eps}=-\gamma(X_t^\eps)\dot{X_t^\eps}+\dot{W_t},
\end{align}
which describes the motion of a particle of mass $\eps$ with the
friction proportional to the velocity. Put $Y_t^\eps=\sqrt{\eps}\dot{X_t^\eps}.$ Then equation (\ref{lav}) can be written as the first order system
\begin{equation} \label{spde15}
\left\{ \begin{aligned}
&\dif X^{\eps}_t =\eps^{-1/2}Y^{\eps}_t\dif t,\\
&\dif Y^{\eps}_t =\eps^{-1}\gamma(X^{\eps}_t)Y^{\eps}_t\dif t+\eps^{-1/2} \dif W_t,
\end{aligned} \right.
\end{equation}
which corresponds to (\ref{spde11}) with $A=F=\Sigma\equiv0$ , $B(x,y)=y$ and $G(x,y)=-\gamma(x)y$. Studying the zero-mass limit  behavior of system (\ref{spde15}) is called the Smoluchovski-Kramers approximation and has been carried out by many authors, see e.g. \cite{CF2,HVW,HMVW,HS} and the references therein.

 \vspace{1mm}
In the finite dimensional situation, the asymptotic behavior for SDEs of the form (\ref{spde11}) was first studied  by Papanicolaou, Stroock and Varadhan \cite{PSV} for a compact state space, see also \cite{Ba} for a similar result in terms of PDEs.  It was found  that  the limit of the slow component will be obtained
in terms of the solution of an auxiliary Poisson equation. Such result  is known as the averaging principle of functional  central limit type, which is also called  the diffusion approximation. Later on, a  non-compact   case  was studied in a series of papers by Pardoux and Veretennikov \cite{P-V,P-V2,P-V3} by using the method of martingale problem, see also  \cite{KY1,RX,RX1} for further development.  To the best of our knowledge, the infinite dimensional setting (\ref{spde11}) has not been studied before.

\vspace{1mm}
To characterize the limit behavior for SPDE (\ref{spde11}), we need to consider the following Poisson equation in the  Hilbert space $H$:
\begin{align}\label{poF}
\cL_2(x,y)\Psi(x,y)=-B(x,y),
\end{align}
where  $\cL_2(x,y)$ is the infinitesimal generator of the frozen process $Y_t^x$ given by (\ref{froz}), i.e.,
\begin{align}\label{L2}
\cL_2\varphi(x,y):=\cL_2(x,y)\varphi(x,y)&:=\langle Ay+G(x,y), D_y\varphi(x,y)\rangle
+\frac{1}{2}Tr\big[D^2_{y}\varphi(x,y)Q_2\big].
\end{align}
It is known that there exists a unique solution $\Psi$ to equation (\ref{poF}) (see Theorem \ref{PP} below).
We shall prove that  the slow process  $X_t^\eps$ in SPDE (\ref{spde11})  converges weakly to $\bar X_t$ as $\eps \to 0$ with $\bar X_t$  solving the following  equation:
\begin{align}\label{spdez}
\dif \bar X_t=A\bar X_t\dif t&+\bar F(\bar X_t)\dif t+\bar\Sigma(\bar X_t)\dif  W_t^1\no\\
&+\overline{B\cdot\nabla_x\Psi}(\bar X_t)\dif t+\Upsilon(\bar X_t)\dif \tilde W_t,
\end{align}
where   $\tilde W_t$ is an $H$-valued cylindrical  Wiener process which is independent of $W^1_t$,
 the new drift coefficient $\overline{B\cdot\nabla_x\Psi}$ and the averaged diffusion coefficient $\bar \Sigma$ are given by
\begin{align*}
\overline{B\cdot\nabla_x\Psi}(x):=\int_{H}\nabla_x\Psi(x,y).B(x,y)\mu^x(\dif y)
\end{align*}
and
 \begin{align}\label{df2}
\<\bar{\Sigma}^2(x)h,k\>:=\int_{H}\<\Sigma(x,y)h,\Sigma(x,y)k\>\mu^x(\dif y),\quad \forall h,k\in H,
\end{align}
and the extra diffusion coefficient $\Upsilon$ is a Hilbert-Schmidt operator satisfying
\begin{align}\label{sst}
\frac{1}{2}\Upsilon(x)\Upsilon^*(x)=\overline{B\otimes\Psi}(x):=\int_{H}\big[B(x,y)\otimes\Psi(x,y)\big]\mu^x(\dif y).
\end{align}
Compared with the averaged equation (\ref{lime}) for SPDE (\ref{spde00}), extra drift term $\overline{B\cdot\nabla_x\Psi}(\bar X_t)\dif t$ and diffusion part $\Upsilon(\bar X_t)\dif \tilde W_t$  appear in (\ref{spdez}), which reflect the homogenization behavior for the fast term $\eps^{-1/2}B(X_t^\eps,Y_t^\eps)\dif t$ in SPDE (\ref{spde11}). Furthermore, we assume that the coefficients in SPDE (\ref{spde11}) are only H\"oler continuous with respect to the fast variable, and we obtain the  rate of convergence of $X_t^\eps$ to $\bar X_t$. Moreover, we deal with the Nemytskii type diffusion and drift coefficients, which require bounds depending on $L^q$-norms and not only on $L^2$-norms. Our result is new even in the case that $\Sigma\equiv0$, and extends the existing results in the literature even in the case $B\equiv0$,  see Remark \ref{br} below for more detailed explanations.

\vspace{1mm}
Our main argument to prove the above convergence is based on the Poisson equation and the Kolmogorov equation in Hilbert space. Undoubtedly, the SPDE (\ref{spde11}) is more difficult than SPDE (\ref{spde00}) due to the presence of the fast term in the slow  equation. Meanwhile, the infinite dimensional situation has more difficulties than the finite dimensional setting, especially in the multiplicative noise case. Some new techniques and nontrivial analysis are needed. First of all,  unlike previous works \cite{Br1,Ce,C2,CF,DS,FWL0,GP1,GP2,GP3,LRS1,LW,PXG,WR}, the uniform moment estimates for $A^\gamma X_t^\eps$ with $\gamma\in[0,1]$  is far from being obvious due to the existence of the fast term $\eps^{-1/2}B(X_t^\eps,Y_t^\eps)\dif t$ in SPDE (\ref{spde11}). In fact, we can only obtain uniform estimates for $A^\gamma X_t^\eps$ with $\gamma\in [0,1/2)$ (which seems to be the best of possible), and the estimates for $A^\gamma X_t^\eps$ with $\gamma\geq 1/2$ will blow-up as $\eps\to0$, see Lemma \ref{la41} below.
Secondly, we need to study  the regularities of the solution of the following infinite dimensional Kolmogorov equation with nonlinear diffusion coefficient:
$$
\partial_t\bar u(t,x)=\bar \cL\,\bar u(t,x),\quad t\in(0,T],
$$
where $\bar\cL$ is the infinitesimal generator of the limit process $\bar X_t$ given by (\ref{spdez}). We point out that even if the diffusion coefficient $\Sigma\equiv0$ in SPDE (\ref{spde11}), we still need to handel the  Kolmogorov equation with nonlinear diffusion coefficient due to the newly generated  diffusion part in SPDE (\ref{spdez}). We also mention that the central limit theorem for SPDE (\ref{spde00}) has been studied in \cite{Ce2,WR} by the classical  time discretization method and in \cite{RXY} by using the Poisson equation, but the  limit processes obtained therein are given by the solutions of  linear equations, which is essentially used in the proof of \cite{RXY}.
Here, we shall need to control terms of the form
 $$
 \<\nabla_x\bar u(t,x),Ax\>\quad\text{and}\quad \<\nabla^2_x\bar u(t,x).Ax,y\>,
 $$
and with $X_t^\eps$ plugged in at the $x$-variable. Even though some new regularities for the infinite dimensional Kolmogorov equations with nonlinear diffusion coefficients have been obtained very recently in \cite{Br4}, the results therein apply only for
 $$
 \<\nabla_x\bar u(t,x),A^\gamma x\>\quad\text{and}\quad \<\nabla^2_x\bar u(t,x).A^\beta x,y\>
 $$
 with $\gamma\in[0,1)$ and $\beta\in[0,1/2)$, which are not sufficient for our purpose. Furthermore, as mentioned above, we do not have uniform control for $A^\gamma X_t^\eps$ with $\gamma\geq 1/2$. For these reasons, we shall use some transfer arguments to handle the low regularities of solution of  the infinite dimensional Kolmogorov equation and the low-order moment estimates for the solution $X_t^\eps$.

\vspace{1mm}
The rest of this paper is organized as follows. In Section 2, we introduce some assumptions and state our main results.  Some preliminaries and uniform estimates for SPDE (\ref{spde11}) are given in Section 3. In Section 4 we give the proof of  the main result.
Throughout this paper, the letter $C$ with or without subscripts
will denote a positive constant, whose value may change in different places, and whose
dependence on parameters can be traced from the calculations.

\vspace{1mm}
{\bf Notations:}
To end this section, we introduce some notations, which will be used throughout this paper. For any $p\in [1,\infty]$,
let $L^p:=L^p(D)$ be the Banach space with  $L^p$-norm $\|\cdot\|_{L^p}$.  In the case of $p = 2$, we denote by $H$ the Hilbert space
$L^2(D)$ endowed with scalar product $\<\cdot,\cdot\>$ and norm $\|\cdot\|.$
For any $p,q\in[2,\infty)$, we use $\sL(L^p,L^q)$ to denote the space of all bounded linear operators from $L^p$ to $L^q.$

\vspace{1mm}
Let $(\gamma_n)_{n\in\mN}$ be a sequence of independent standard Gaussian
random variables.
An operator $\Phi\in\sL(H,L^q)$ is said to be  $\gamma$-Radonifying (see e.g. \cite[Section 2.1]{Br2}) if there exists an orthonormal system $(e_n)_{n\in \mN}$ of $H$ such that the series $\sum\limits_{n\in\mN}\gamma_n\Phi e_n$ converges in $L^2(\Omega,L^p).$
 We shall denote by $\sR(H,L^q)$ the space of all $\gamma$-Radonifying operators from $H$ to $L^p,$ with the norm $\|\cdot\|_{\sR(H,L^q)}$ defined by
$$\|\Phi\|_{\sR(H,L^q)}:=\mE\Big\|\sum\limits_{n\in\mN}\gamma_n\Phi e_n\Big\|^2_{L^p}.$$
For any $p\in[2,\infty),\Phi\in\sR(H,L^p),$ it is known that there exists a constant $C_p>0$ such that
\begin{align}\label{prot1}\|\Phi\|_{\sR(H,L^p)}\leq C_p\Big\|\sum\limits_{n\in\mN}(\Phi e_n)^2\Big\|_{L^{p/2}}.
\end{align}
When $p=2$, $\sR(H,H)=\sL_2(H)$  is the space of all Hilbert-Schmidt operators on $H$ and $\|\Phi\|_{\sR(H,H)}=Tr(\Phi\Phi^*).$
Let $(W_t)_{t\geq0}$ be an $H$-valued  Wiener process. Then, for any $T\in[0,\infty),p\in[2,\infty)$  and predictable processes $\Phi\in L^2(\Omega\times[0,T];\sR(H,L^p))$, the $L^p$-valued It\^o integral
$\int_0^T\Phi(t)\dif W_t$
is well defined. Moreover, there exists $C_p>0$ such that
\begin{align}\label{prot2}
\mE\bigg(\bigg\|\int_0^T\Phi(t)\dif W_t\bigg\|^2_{L^p}\bigg)\leq C_p\bigg(\int_0^T\mE\|\Phi(t)\|^2_{\sR(H,L^p)}\dif t\bigg).
\end{align}

\vspace{1mm}
For any $x,y\in H$ and $\phi: H\times H\to \hat H$, where $\hat H$ is another Hilbert space,  we say that $\phi$ is G\^ateaux differentiable at $x$ if there exists an operator $D_x\phi(x,y)\in \sL(H,\hat H)$ such that for all $h\in H$,
$$
\lim_{\tau\to 0}\frac{\phi(x+\tau h,y)-\phi(x,y)}{\tau}=D_x\phi(x,y).h.
$$
If in addition
$$
\lim_{\|h\|\to 0}\frac{\|\phi(x+h,y)-\phi(x,y)-D_x\phi(x,y).h\|_{\hat H}}{\|h\|}=0,
$$
$\phi$ is called Fr\'echet differentiable at $x$. Similarly, for any $k\geq 2$ we can define the $k$ times G\^ateaux and Fr\'echet derivative of $\phi$ at $x$, and we will identify the higher order  derivatives $D^k_x\phi(x,y)$ with a linear operator in $\sL^k(H,\hat H):=\sL(H,\sL^{(k-1)}(H,\hat H))$, endowed with the operator norm
\begin{align*}
\|D_x^k\phi(x,y)\|_{\sL^k(H,\hat H)}:=\sup_{\|h_1\|\leq 1,\|h_2\|\leq1,\cdots,\|h_k\|\leq1, \|h\|\leq 1}\<D^k_x\phi(x,y).(h_1,h_2,\cdots,h_k),h\>_{\hat H}.
\end{align*}
By the same way, we define the G\^ateaux and Fr\'echet derivatives of $\phi$ with respect to the $y$ variable, and we have $D_y\phi(x,y)\in \sL(H,\hat H)$, and for $k\geq 2$, $D^k_y\phi(x,y)\in \sL^k(H,\hat H) :=\sL(H,\sL^{(k-1)}(H,\hat H))$.

\vspace{1mm}
We will denote by $L^\infty(H\times H,\hat H)$  the space of  all measurable maps $\phi: H\times H\to \hat H$ satisfying
$$
\|\phi\|_{L^\infty(\hat H)}:=\sup_{(x,y)\in H\times H}\|\phi(x,y)\|_{\hat H}<\infty.
$$
For $k\in\mN$, the space $C_b^{k,0}(H\times H,\hat H)$ consists of  all maps $\phi\in L^\infty(H\times H,\hat H)$ which are $k$ times G\^ateaux differentiable at any $x\in H$ with bounded derivatives.
Similarly, the space $C_b^{0,k}(H\times H,\hat H)$ consists of all maps $\phi\in L^\infty(H\times H,\hat H)$  which are $k$ times G\^ateaux differentiable at any $y\in H$  with bounded derivatives.
 We also introduce the space $\mC_b^{0,k}(H\times H,\hat H)$ consisting of all maps which are $k$ times Fr\'echet differentiable at any $y\in H$ with bounded derivatives.

\vspace{1mm}
For $\eta\in(0,1)$, we use $C_b^{k,\eta}(H\times H,\hat H)$ to denote the subspace of $C_b^{k,0}(H\times H,\hat H)$ consisting of all maps such that
$$
\|\phi(x,y_1)-\phi(x,y_2)\|_{\hat H}\leq C_0\|y_1-y_2\|^\eta.
$$
When $\hat H=\mR$, we will omit the letter $\hat H$ in the above notations  for simplicity.

\section{Assumptions and Main results}

Let $\{e_{n}\}_{n\in \mN}$ be a complete orthonormal basis of $H$. Throughout this paper, we assume that there exist non-decreasing sequences of real positive numbers $\{\alpha_n\}_{n\in\mN}$ such that
\begin{align}\label{aann}
Ae_{n}=-\alpha_ne_{n},\quad \forall n\in\mN.
\end{align}
In this setting, the power of  $-A$ can be easily defined as follows: for any $\theta\in[0,1]$,
$$
(-A)^\theta x:=\sum_{n\in\mN}\alpha_n^{\theta}\langle x,e_{n}\rangle e_{n},
$$
with the domain
$$
\cD((-A)^\theta):=\bigg\{x\in H: \|x\|_{(-A)^\theta}^2:=\sum\limits_{n\in\mN}\alpha_n^{2\theta}\langle x,e_{n}\rangle^2<\infty\bigg\}.
$$
Moreover, the corresponding semigroup $\{e^{tA}\}_{t\geq0}$ can be defined through the following spectral formula: for any $t\geq 0$ and $x\in H$,
$$
e^{tA}x:=\sum_{n\in\mN}e^{-\alpha_nt}\langle x,e_{n}\rangle\,e_{n}.
$$
Then it is known that for any $\gamma\in[0,1],t>0$ and $p\in[2,\infty)$, we have (see e.g. \cite[(3)]{Br2})
\begin{align}\label{pp1}
\|(-A)^\gamma e^{tA}x\|_{L^p}\leq C_{\gamma,p}t^{-\gamma}e^{-\frac{\alpha_1}{2}t}\| x\|_{L^p},
\end{align}
where $C_{\gamma,p}>0$ is a constant. Furthermore, for any $\theta\in[0,1/4),$  it follows from \cite[(10)]{Br4} (see also \cite[Proposition 2.1]{Br2}) that for any Lipschitz continuous function $\varphi$ and $x\in\cD_p((-A)^{\theta+\delta})$ with $\delta>0,$
\begin{align}\label{lfp}
\|(-A)^\theta\varphi(x)\|_{L^p}\leq C_{\theta,\delta,\varphi}\big(1+\|(-A)^{\theta+\delta}x\|_{L^p}\big).
\end{align}

For $i=1,2$, let $Q_i$ be two linear self-adjoint bounded operators on $H$ with positive eigenvalues $\{\lambda_{i,n}\}_{n\in\mN}$, i.e.,
$$
Q_ie_{n}=\lambda_{i,n}e_{n},\quad\forall n\in\mN.
$$
 Recall that $W_t^i,i=1,2$, are $H$-valued $Q_i$-Wiener processes both defined on a complete filtered probability space $(\Omega,\mathscr{F},\mathscr{F}_t,\mP)$. It is known that $W_t^i$ can be written as
$$
W_t^i=\sum_{n\in\mN}\sqrt{\lambda_{i,n}}\beta_{i,n}(t)e_{i,n},
$$
where $\{\beta_{i,n}\}_{n\in\mN}$ are mutual independent real-valued Brownian motions. We shall always assume that for $i=1,2,$ $Q_i$ are trace operators, and for $\gamma\in[0,1/2)$ and $p\in[2,\infty)$,
\begin{align}\label{asuq}
\int_0^T\big\|(-A)^\gamma e^{tA}Q_i^{1/2}\big\|^2_{\sR(L^2,L^p)}\dif t<\infty,
\end{align}
and for any $T>0$ , we have
\begin{align}
\int_0^T\Lambda_t^{\frac{1+\vartheta}{2}}\dif t<\infty,
\end{align}
where
\begin{align}
\Lambda_t:=\sup\limits_{n\geq 1}\frac{2\alpha_n}{\lambda_{2,n}(e^{2\alpha_nt}-1)}<\infty,
\end{align}
$\alpha_n$ is given by (\ref{aann}), and $\vartheta\geq\max{(\eta,1-\eta)}$ with $\eta$ being the H\"older regularity of the coefficients in the assumption of Theorem 2.1 below.

Furthermore, we assume that $B: H\times H\to H$ and $\Sigma: H\times H\to\sL(H)$ are defined as the Nemytskii operators, i.e., there exist $b, \sigma: \mR\times \mR\to\mR$  such that
\begin{align}\label{BS}
B(x,y)(\xi)=b(x(\xi),y(\xi))\quad\text{and}\quad[\Sigma(x,y)z](\xi)=\sigma(x(\xi),y(\xi))z(\xi).
\end{align}
We also assume that
$B$ satisfies the centering condition:
\begin{align}\label{mcen}
\int_H B(x,y)\mu^x(\dif y)=0,
\end{align}
where $\mu^x(\dif y)$ is the unique invariant measure  of the  frozen process $Y_t^x$. Such kind of assumption is necessary and analogous to the centering in the standard central limit theorem, see e.g. \cite{PSV,P-V,P-V2,P-V3}.

The following is the main result of this paper.

\bt\label{main3}
 Let $T>0$ and $x,y\in L^8$. Assume that  $f,b,g,\sigma\in C^{4,\eta}_b(\mR^2,\mR)$ with $\eta>0$. Then for any $\varphi\in C_b^4(H)$,   we have
\begin{align}\label{555}
\sup_{t\in[0,T]}\big|\mE[\varphi(X_t^{\eps})]-\mE[\varphi(\bar X_{t})]\big|\leq C_0\,\eps^{\frac{1}{2}},
\end{align}
where  $C_0=C(T,x,y,\varphi)>0$ is a constant independent of $\eta$ and $\eps$.
\et

We list some important comments to explain our result.

\br\label{br}

(i) Our result seems to be new even when $B\equiv0$. In fact, as far as we know, the  multiplicative noise case of SPDE (\ref{spde11}) where the diffusion coefficient depends on both the fast and the slow variables has been studied only in \cite{Ce} when $B\equiv0$. The argument  in \cite{Ce} is based on the classical Khasminskii's time discretisation approach and no  rate of convergence is obtained therein.
In the present paper, by following exactly the same procedure as in our proof (in fact, more easily if $B\equiv0$), we can get that
\begin{align}\label{777}
\sup_{t\in[0,T]}\big|\mE[\varphi(X_t^{\eps})]-\mE[\varphi(\bar X_{t})]\big|\leq C_0\,\eps.
\end{align}
This means that the rate is of order 1 in the weak convergence of the averaging principle, which coincides with the finite dimensional situation, see e.g. \cite{Br1,KY1}.

\vspace{1mm}
\noindent
(ii)  We point out that  the noise part in the slow equation can be totally degenerate, i.e., we allow $\Sigma\equiv0$ in SPDE (\ref{spde11}). Even in this  case, the limit behavior  for equation (\ref{spde11}) has not been studied before in the infinite dimensional situation due to the  appearance of the fast term. Unlike the  convergence in the averaging principle of SPDE (\ref{spde11}) with $B\equiv0$, where the noise in the limit equation is  additive     if the original slow equation is driven by additive noise  (see e.g. \cite{Br1,C2,CF,DS,FWL0,GP1,GP2,GP3,LRS1,LW,PXG,WR}), the main difference now is that even though the noise  is additive or $\Sigma\equiv0$ (totally degenerate) in SPDE (\ref{spde11}), the corresponding limit equation will exhibit multiplicative noise  in view of the newly generated diffusion part in (\ref{spdez}). This is due to the  homogenization effect of the fast term in SPDE (\ref{spde11}).

\vspace{1mm}
\noindent
(iii)  The 1/2 order rate of convergence in (\ref{555}) is known to be optimal in the finite dimensional situation in view of the asymptotic expansion in \cite{KY1}. Intuitively, the difference between (\ref{777}) and (\ref{555}) is caused by the fast term $\frac{1}{\sqrt{\eps}}B$, which reduces the convergence rate from $\eps$ to $\sqrt{\eps}$.
Note that we assume   the coefficients of SPDE (\ref{spde11}) are only $\eta$-H\"older with respect to the fast variable, and  the rate of convergence does not depend on $\eta$. This reflects that the slow process is the main term in the limiting procedure of the multi-scale system, which coincides  with  intuition since the fast component has been totally homogenized out in the limit equation.
\er

\section{Preliminaries and a priori estimates}

In this section, we prove some uniform estimates, with respect to $\eps\in(0,1)$, for the solution $(X_t^\eps, Y_t^\eps)$ of system (\ref{spde11}). In fact, the  estimates for the fast variable $Y_t^\eps$ can be proved similarly as in previous works. However, the uniform control for the slow variable $X_t^\eps$ is far from being obvious due to the existence of the fast term $\eps^{-1/2}B(X_t^\eps,Y_t^\eps)$ in the  equation. For this, we shall derive some strong fluctuation estimates by using the technique of  Poisson equation.

\subsection{Preliminaries}

Recall that the drift coefficients $F, G$ and $B$ are  Nemytskii operators defined by (\ref{FG}) and (\ref{BS}), respectively, and we assumed that $f,b,g\in C^{4,\eta}_b(\mR^2,\mR)$ with $\eta>0$. However, it is well-known that $F, G$ and  $B$ do not inherit higher order  regularity properties on $H$. The control of their higher order derivatives requires the use of $L^p$ norms. The following properties can be found in \cite[Property 3.2]{Br4}. We write them for $F,$ but they also hold with $G$ and $B.$

\bl\label{fbg}
Let $F(\cdot,\cdot): H\times H\to H$ be defined by (\ref{FG}). Then for every $y\in H$, $F(\cdot,y)$ is  fourth times G\^ateaux differentiable. Moreover, the following properties hold:

\vspace{1mm}
	\noindent (i) $F\in C_b^{1,\eta}(H\times H, H);$

\vspace{1mm}
\noindent (ii) for any $x,y\in H$ and $p,r_1,r_2\in[1,\infty]$ satisfying $\frac{1}{p}=\frac{1}{r_1}+\frac{1}{r_2},$
$$\|D_x^2F(x,y).(h_1,h_2)\|_{L^p}\leq C_1\,\|h_1\|_{L^{r_1}}\|h_2\|_{L^{r_2}};$$

\vspace{1mm}
\noindent (iii) for any $x,y\in H$ and $p,q_1,q_2,q_3\in[1,\infty]$ satisfying $\frac{1}{p}=\frac{1}{q_1}+\frac{1}{q_2}+\frac{1}{q_3},$
$$\|D_x^3F(x,y).(h_1,h_2,h_3)\|_{L^p}\leq C_2\,\|h_1\|_{L^{q_1}}\|h_2\|_{L^{q_2}}\|h_3\|_{L^{q_3}};$$

\vspace{1mm}
\noindent (iv) for any $x,y\in H$ and $p,q_1,q_2,q_3,q_4\in[1,\infty]$ satisfying $\frac{1}{p}=\frac{1}{q_1}+\frac{1}{q_2}+\frac{1}{q_3}++\frac{1}{q_4},$
$$\|D_x^4F(x,y).(h_1,h_2,h_3,h_4)\|_{L^p}\leq C_3\,\|h_1\|_{L^{q_1}}\|h_2\|_{L^{q_2}}\|h_3\|_{L^{q_3}}\|h_4\|_{L^{q_4}},$$
where $C_i,i=1,2,3$ are positive constants.
\el

\vspace{1mm}
As for the diffusion coefficient $\Sigma$ defined by (\ref{BS}), due to $\sigma\in C_b^{4,\eta}(\mR^2,\mR)$ with $\eta>0$,  we have the following result, see e.g. \cite[Property 3.3]{Br4}.

\bl\label{sig}
Let $\Sigma(\cdot,\cdot):H\times H\to\sL(H)$ be defined by (\ref{BS}). Then for every $y\in H$, $\Sigma(\cdot,y)$ is fourth times G\^ateaux differentiable. Moreover, we have:

\vspace{1mm}
\noindent (i) for any $x,y\in H,\|\Sigma(x,y)\|_{\sL(H)}\leq C_1;$

\vspace{1mm}
\noindent (ii) for any $x,y\in H$ and $h\in L^\infty,$
$\|D_x\Sigma(x,y).h\|_{\sL(H)}\leq C_2\,\|h\|_\infty;$

\vspace{1mm}
\noindent (iii) for $x,y\in H$ and $h_1,h_2\in L^\infty,$
$\|D_x^2\Sigma(x,y).(h_1,h_2)\|_{\sL(H)}\leq C_3\,\|h_1\|_\infty\|h_2\|_\infty;$

\vspace{1mm}
\noindent (iv) for $x,y_1,y_2\in H$ and $h\in L^\infty,$
$\|[\Sigma(x,y_1)-\Sigma(x,y_2)].h\|_{\sL(H)}\leq C_4\,\|y_1-y_2\|^\eta\|h\|_\infty,$
where $C_i,i=1,2,3,4$ are positive constants.
\el

Consider the following Poisson equation in the infinite dimensional Hilbert space $H$:
\begin{align}\label{pois}
\cL_2(x,y)\psi(x,y)=-\phi(x,y),
\end{align}
where $\cL_2(x,y)$ is defined by (\ref{L2}), $x\in H$ is regarded as a parameter, and $\phi: H\times H\rightarrow \hat H$ is measurable.
Recall that $Y_t^x(y)$  satisfies the frozen equation (\ref{froz}) and $\mu^x(\dif y)$ is the (unique)  invariant measure of  $Y_t^x(y)$.
Since there is no boundary condition in (\ref{pois}), to be well-posed, we need to make the following ``centering" assumption  on $\phi$:
\begin{align}\label{cen}
\int_{H}\phi(x,y)\mu^x(\dif y)=0,\quad\forall x\in H.
\end{align}
The following result has been proven in \cite[Theorem 3.2]{RXY}.

\bt\label{PP}
For every $\phi:H\times H\to\hat H$ satisfying  (i)-(ii) of Lemma \ref{fbg} and the centering condition (\ref{cen}),
there exists  a  unique classical solution to the equation (\ref{pois})  which is given by
\begin{align*}
\psi(x,y)=\int_0^\infty\!\mE\big[\phi(x,Y_t^x(y))\big]\dif t,
\end{align*}
where $Y_t^x(y)$ satisfies the frozen equation (\ref{froz}). Moreover, we have

\vspace{1mm}
\noindent
(i) $\psi\in \mC_b^{0,2}(H\times H,\hat H);$

\vspace{1mm}
\noindent
(ii) $\psi$ is twice G\^ateaux differentiable with respect to the first variable, and the derivatives satisfy estimates in (i)-(ii) of Lemma \ref{fbg}.
\et
\br
According to   \cite[Lemma 3.7]{RXY},  we also have that $\bar F\in C_b^1(H,H)$ and the $k$-th (k=2,3,4) derivatives  satisfy estimates in (ii)-(iv) of Lemma \ref{fbg}.
\er
We shall need to use  It\^o's formula for $\psi(x,y)$ with $(X_t^\eps, Y_t^\eps)$ plugged in at both variables, say $\psi(X_t^\eps, Y_t^\eps)$.  However, due to the presence of the unbounded operator in equation (\ref{spde11}) and the fact that $\psi$ is only G\^ateaux differentiable with respect to the $x$-variable, we can not apply It\^o's formula for SPDE (\ref{spde11}) directly. For this reason, we recall the following Galerkin approximation scheme.

For $n\in\mN$, let $H^n:= span\{e_{k};1\leq k\leq n\}$ and denote  the orthogonal projection of $H$ onto $H^n$ by $P_n$. For $(x,y)\in H^n\times H^n$, define
$$F_n(x,y):=P_nF(x,y),\;B_n(x,y):=P_nB(x,y),\;G_n(x,y):=P_nG(x,y).$$
We reduce the infinite dimensional system (\ref{spde11}) to the following finite dimensional system in $H^n\times H^n$:
\begin{equation}\label{xyz}
\left\{ \begin{aligned}
&\dif X^{n,\eps}_t =AX^{n,\eps}_t\dif t+F_n(X^{n,\eps}_t, Y^{n,\eps}_t)\dif t\\
&\qquad\qquad\qquad\quad\,\,\,+\eps^{-1/2}B_n(X^{n,\eps}_t, Y^{n,\eps}_t)\dif t+P_n\Sigma(X^{n,\eps}_t, Y^{n,\eps}_t)\dif W_t^1,\\
&\dif Y^{n,\eps}_t =\eps^{-1}AY^{n,\eps}_t\dif t+\eps^{-1}G_n(X^{n,\eps}_t, Y^{n,\eps}_t)\dif t+\eps^{-1/2} P_n\dif W_t^2,
\end{aligned} \right.
\end{equation}
with initial values $X_0^{n,\eps}=x^n:=P^nx\in H^n$ and $Y_0^{n,\eps}=y^n:=P_ny\in H^n$.
It is easy to check that $F_n,B_n$ and $G_n$ satisfy the same conditions as $F,B$ and $G$ with bounds which are uniform with respect to $n$. The corresponding averaged equation for system (\ref{xyz}) can be formulated as
\begin{align*}
\dif \bar{X}^n_t=A\bar{X}^n_t\dif t+\bar{F}_n(\bar{X}^n_t)\dif t+\overline{(B\cdot\nabla_x\Psi)}_n(\bar X^n_t)\dif t+P_n\Upsilon(\bar X^n_t)\dif \tilde W_t&+P_n\bar\Sigma(\bar X^n_t)\dif  W_t^1,
\end{align*}
where $\bar{F}_n(x),\overline{(B\cdot\nabla_x\Psi)}_n(x) $ are defined by
\begin{align*}
\bar F_n(x):=\int_{H^n}F_n(x,y)\mu^x_n(\dif y)
\end{align*}
and
\begin{align*}
\overline{(B\cdot\nabla_x\Psi)}_n(x)=:\int_{H^n}D_x\Psi_n(x,y).B_n(x,y)\mu^x_n(\dif y),
\end{align*}
respectively. For any $T>0$ and $\varphi\in C_b^4(H),$ we have for $t\in[0,T]$,
\begin{align}\label{wcz}
\left|\mE[\varphi(X_t^{\eps})]-\mE[\varphi(\bar X_t)]\right|&\leq \left|\mE[\varphi(X_t^{\eps})]-\mE[\varphi(X_t^{n,\eps})]\right|\no\\
&\quad+\left|\mE[\varphi(X_t^{n,\eps})]-\mE[\varphi(\bar X_t^n)]\right|+\left|\mE[\varphi(\bar X_t^n)]-\mE[\varphi(\bar X_t)]\right|.
\end{align}
By using similar arguments as in  \cite[Lemma 5.4]{RXY}, the first and the last terms on the right-hand of (\ref{wcz}) converge to $0$ as $n\to\infty$. Therefore,
in order to prove Theorem {\ref{main3}}, we only need to show that
\begin{align}\label{nzz}
\sup_{t\in[0,T]}\left|\mE[\varphi(X_t^{n,\eps})]-\mE[\varphi(\bar X_t^n)]\right|\leq C_T\,\eps^{\frac{1}{2}},
\end{align}
where {\bf $C_T>0$ is a constant independent of} $n.$ In the rest part of this paper, we shall only work with the approximating system (\ref{xyz}), and proceed to prove
bounds that are uniform with respect to $n$. To simplify the notations, we shall omit the index $n$. In particular, the space $H^n$ are denoted by $H.$

\subsection{Moment estimates}
Let $(X_t^\eps,Y_t^\eps)$ satisfy the following equation:
\begin{equation} \label{spde22}
	\left\{ \begin{aligned}
	&X^{\eps}_t =e^{tA}x+\int_0^te^{(t-s)A}F(X^{\eps}_s, Y^{\eps}_s)\dif s+\eps^{-1/2}\int_0^te^{(t-s)A}B(X^{\eps}_s, Y^{\eps}_s)\dif s\\&\qquad\qquad\qquad\qquad\qquad\qquad\qquad\qquad+\int_0^te^{(t-s)A}\Sigma(X_s^\eps,Y_s^\eps)\dif W_s^1,\\
	& Y^{\eps}_t =e^{\frac{t}{\eps}A}y+\eps^{-1}\int_0^te^{\frac{t-s}{\eps}A}G(X^{\eps}_s, Y^{\eps}_s)\dif s+\eps^{-1/2}\int_0^te^{\frac{t-s}{\eps}A}\dif W_s^2.\\
	\end{aligned} \right.
	\end{equation}
Recall that $\cL_2(x,y)$ is defined by (\ref{L2}). For convenience, we  denote by
\begin{align}\label{L}
\cL\varphi(x,y):=\cL_1\varphi(x,y)+\eps^{-1/2}\cL_0\varphi(x,y),\quad\forall \varphi\in C_b^{2,0}(H\times H),
\end{align}
where
\begin{align}\label{L1}
\cL_1\varphi(x,y):=\cL_1(x,y)\varphi(x,y):=&\langle Ax+F(x,y), D_x\varphi(x,y)\rangle\no\\
&
+\frac{1}{2}Tr\big[D_x^2\varphi(x,y)\Sigma(x,y)Q_1\Sigma^*(x,y)\big],
\end{align}
and
\begin{align}\label{L0}
\cL_0\varphi(x,y):=\cL_0(x,y)\varphi(x,y):=\<B(x,y),
D_x\varphi(x,y)\>.
\end{align}
The following moment estimates  for the fast variable $Y_t^\eps$  can be proved by using the similar arguments as in \cite[Propositions A.2 and A.4]{Br1} and the properties (\ref{prot1}) and (\ref{prot2}), we omit the details here. 

\begin{lemma}\label{la42}
Let $T>0$ and  $y\in L^p$ with $p\in[2,\infty)$.  Then

\vspace{1mm}
	\noindent (i) for any $q\geq1$, $\gamma\in[0,1/2)$ and $t\in(0,T]$, we have
\begin{align}\label{msy}
	\sup\limits_{\eps\in(0,1)}\mE\|(-A)^\gamma Y^{\eps}_t\|_{L^p}^{q}\leq C_{\gamma,p,q,T}\,t^{-\gamma q}\Big(1+\| y\|_{L^p}^{q}\Big);
	\end{align}
(ii)  for  any  $q\geq1$, $\gamma\in[0,1/2]$ and $0< s\leq t\leq T,$	 we have
\begin{align}\label{msy2}
\big(\mE\|  Y_t^\eps- Y_s^\eps\|_{L^p}^{q}\big)^{\frac{1}{q}}\leq C_{\gamma,p,q,T}\,\bigg(\frac{(t-s)^{\gamma}}{s^{\gamma}}e^{-\frac{\alpha_1}{2\eps}s}
\|y\|_{L^p}
+\frac{(t-s)^{\gamma}}{\eps^{\gamma}}\bigg);
	\end{align}
where $C_{\gamma,p,q,T}>0$ is a constant.
\end{lemma}

Concerning the estimates for $X_t^\eps,$ by regarding the term $F+\eps^{-1/2}B$ as the whole drift coefficient and following  exactly the same arguments as in \cite[Proposition 2.10]{Br2}, we easily have the following preliminary results.

\bl\label{laxts}Let $T>0$ and  $x\in L^p$ with $p\in[2,\infty)$.  Then

\vspace{1mm}
	\noindent (i)
for any $q\geq1$, $\gamma\in[0,1/2)$ and $t\in(0,T]$, we have
\begin{align}\label{msx}
	\sup\limits_{\eps\in(0,1)}\mE\| (-A)^\gamma X^{\eps}_t\|_{L^p}^{q}\leq C_{\gamma,p,q,T}\,t^{-\gamma q}\,\eps^{-q/2}\big(1+\| x\|_{L^p}^{q}\big);
	\end{align}
(ii) for  any  $q\geq1$, $\gamma\in[0,1/2]$ and $0< s\leq t\leq T,$	 we have
\begin{align}\label{xts}
	\Big(\mE\|  X_t^\eps- X_s^\eps\|_{L^p}^q\Big)^{\frac{1}{q}}&\leq C_{\gamma,p,q,T}\,\bigg(\frac{{(t-s)}^{\gamma}}{s^{\gamma}}e^{-\frac{\alpha_1}{2}s}
\|x\|_{L^p} +\frac{(t-s)^{\gamma}}{\eps^{1/2}}\bigg);
	\end{align}
where $C_{\gamma,p,q,T}>0$ is a constant.
\el

However, the moment estimate (\ref{msx}) is not enough to use below since it blows up as $\eps\to0$. We need some uniform estimates for $X_t^\eps$ with respect to $\eps\in(0,1)$. For this, we establish the following strong fluctuation estimate for the integral functional of $(X_s^\eps,Y_s^\eps)$ over the time interval $[0,t].$

\bl[Strong fluctuation estimate]\label{strf}
Let $T>0$ and $x,y\in L^p$ with $p\in[2,\infty)$. Then for any $\gamma\in[0, 1/2),$ $q\geq 1$, $0\leq t\leq T$ and $\phi :H\times H \to H$ satisfying both (i)-(ii) of Lemma \ref{fbg} and the centering condition (\ref{cen}), we have
\begin{align}\label{stfe}
\mE\left\|\int_0^t(-A)^\gamma e^{(t-s)A}\phi(X_s^\eps,Y_s^\eps)\dif s\right\|_{L^p}^q\leq C_{\gamma,p,q,T}\,t^{-\gamma q}\,\eps^{q/2}(1+\|x\|_{L^p}^{q}+\|y\|_{L^p}^{q}),
\end{align}
where $C_{\gamma,p,q,T}>0$ is a constant.
\el

\begin{proof}
Let $\psi$  solve the Poisson equation
$$
\cL_2(x,y)\psi(x,y)=-\phi(x,y),
$$
and define
\begin{align}\label{psi8}
\psi_{t,\gamma}(s,x,y):=(-A)^\gamma e^{(t-s)A}\psi(x,y).
\end{align}
Note that $\cL_2$ is an operator with respect to the $y$-variable, it is easy to verify that
\begin{align}\label{ppo}
\cL_2(x,y)\psi_{t,\gamma}(s,x,y)=-(-A)^\gamma e^{(t-s)A}\phi(x,y).
\end{align}
 In view of Theorem \ref{PP}, we can apply It\^o's formula to $\psi_{t,\gamma}(t, X_t^\eps,Y_t^{\eps})$ to get
\begin{align}\label{ito1}
\psi_{t,\gamma}(t, X_t^\eps,Y_t^{\eps})&=\psi_{t,\gamma}(0, x,y)+\int_0^t (\p_s+\eps^{-1/2}\cL_0+\cL_1)\psi_{t,\gamma}(s,X_s^\eps,Y_s^{\eps})\dif s\no\\
&\quad+\frac{1}{\eps}\int_0^t\cL_2\psi_{t,\gamma}(s,X_s^\eps,Y_s^{\eps})\dif s+M_{t}^1+\frac{1}{\sqrt{\eps}}M_{t}^2,
\end{align}
where $M_{t}^1$ and $M_{t}^2$ are  defined by
\begin{align*}
M_{t}^1:=\int_0^t \<D_x\psi_{t,\gamma}(s,X_s^\eps,Y_s^{\eps}),\Sigma(X_s^\eps,Y_s^{\eps})\dif W_s^1\>\end{align*}
and
\begin{align*}
 M_{t}^2:=\int_0^tD_y\psi_{t,\gamma}(s,X_s^\eps,Y_s^{\eps})\dif W_s^2.
\end{align*}
Multiplying  both sides of (\ref{ito1}) by $\eps$ and using (\ref{ppo}), we get
\begin{align*}
\int_0^t(-A)^\gamma e^{(t-s)A}\phi(X_s^\eps,Y_s^\eps)\dif s&=-\int_0^t\cL_2\psi_{t,\gamma}(s,X_s^\eps,Y_s^{\eps})\dif s\\
&=\eps\,\big[\psi_{t,\gamma}(0, x,y)-\psi_{t,\gamma}(t,X^\eps_t,Y_t^{\eps})\big]+\eps\, M_{t}^1+\sqrt{\eps}\,M_{t}^2\\
&+\eps\int_0^t(\p_s+\eps^{-1/2}\cL_0+\cL_1)\psi_{t,\gamma}(s,X_s^\eps,Y_s^{\eps})\dif s\\
&=\eps\, (-A)^{\gamma}e^{tA}\big[\psi(x,y)-\psi(X^\eps_t,Y_t^{\eps})\big]\\
&+\eps\int_0^t(-A)^{1+\gamma} e^{(t-s)A}\big[\psi(X^\eps_s,Y_s^{\eps})-\psi(X^\eps_t,Y_t^{\eps})\big]\dif s\\
&+\eps\int_0^t(\eps^{-1/2}\cL_0+\cL_1)\psi_{t,\gamma}(s,X_s^\eps,Y_s^{\eps})\dif s+\eps\, M_{t}^1+\sqrt{\eps}\,M_{t}^2.
\end{align*}
For any $0\leq t \leq T$ and $q\geq 1$,
we deduce that
\begin{align*}
&\mE\left\|\int_0^t(-A)^\gamma e^{(t-s)A}\phi(X_s^\eps,Y_s^\eps)\dif s\right\|_{L^p}^q\\
&\leq C_0\,\bigg(\eps^q\,\mE\big\| (-A)^{\gamma}e^{tA}\big[\psi(x,y)-\psi(X^\eps_t,Y_t^{\eps})\big]\big\|_{L^p}^q\\
&\quad+\eps^q\,\mE\left\|\int_0^t (-A)^{1+\gamma} e^{(t-s)A}\big[\psi(X^\eps_s,Y_s^{\eps})-\psi(X^\eps_t,Y_t^{\eps})\big]\dif s\right\|_{L^p}^q\\
&\quad+\eps^{q/2}\,\mE\left\|\int_0^t \cL_0\psi_{t,\gamma}(s,X_s^\eps,Y_s^{\eps})\dif s\right\|_{L^p}^q\\
&\quad+\eps^q\,\mE\left\|\int_0^t \cL_1\psi_{t,\gamma}(s,X_s^\eps,Y_s^{\eps})\dif s\right\|_{L^p}^q+\eps^q\,\mE\|M_{t}^1\|_{L^p}^q+\eps^{q/2}\,\mE\|M_{t}^2\|_{L^p}^q\bigg)\!=:\!\sum_{i=1}^6\sJ_i(t,\eps).
\end{align*}
For the first term, by (\ref{pp1}) and Theorem \ref{PP},
we have
\begin{align*}
\sJ_{1}(t,\eps)&\leq C_1\,\eps^q\,t^{-\gamma q}.
\end{align*}
For $\gamma'\in(\gamma,1/2)$, it follows from (\ref{msy2}) and (\ref{xts}) that
\begin{align*}
\sJ_{2}(t,\eps)\!
&\leq C_2\,\eps^q\Bigg(\!\!\int_0^t(t-s)^{-1-\gamma}\bigg[\Big(\mE\big[\| X_t^\eps- X_s^\eps\|_{L^p}^{2q}\big]\Big)^{1/2q}\!\!+\!\Big(\mE\big[\|Y_t^{\eps}-Y_s^{\eps}\|_{L^p}^{2 q}\big]
\Big)^{1/2q}\bigg]\dif s\Bigg)^q\\
&\leq C_2\,\eps^q(1+\|x\|_{L^p}^q+\|y\|_{L^p}^q)\left(\int_0^t(t-s)^{-1-\gamma}\frac{(t-s)^{\gamma'}}{\eps^{1/2}}\dif r\right)^q\\
&\leq  C_2\,\eps^{q/2}(1+\|x\|_{L^p}^q+\|y\|_{L^p}^q).
\end{align*}
For the third term,  by  definition (\ref{L0}) and Theorem \ref{PP}, we have
\begin{align*}
\sJ_{3}(t,\eps)
&\leq C_3\;\eps^{q/2}\int_0^t(t-s)^{-\gamma}\dif s\leq C_3\;\eps^{q/2}.
\end{align*}
To deal with the fourth term, by definitions (\ref{L1}), (\ref{psi8}), Minkowski's inequality, Lemma \ref{laxts} and Theorem {\ref{PP}}, we deduce that for $\gamma'\in(\gamma,1/2),$
\begin{align*}
\sJ_{4}(t,\eps)&\leq C_4\,\eps^q\,\left(\int_0^t \big(\mE\|\cL_1\psi_{t,\gamma}(s,X_s^\eps,Y_s^{\eps})\|_{L^p}^q\big)^{1/q}\dif s\right)^q\\&\leq  C_4\,\eps^{q}+ C_4\,\eps^q\,\left(\int_0^t \big(\mE|\< (-A)^{\gamma'}X_s^\eps,(-A)^{1-\gamma'+\gamma}e^{(t-s)A}D_x\psi(X_s^\eps,Y_s^\eps)\>|^q\big)^{1/q}\dif s\right)^q
\\&\leq C_4\,\eps^{q}+ C_4\,\eps^q\left(\int_0^t(t-s)^{-1+\gamma'-\gamma}(\mE\|(-A)^{\gamma'}X_s^\eps\|_{L^p}^{q})^{1/q}\dif s\right)^q\\
&\leq C_4\,\eps^{q/2}(1+\|x\|_{L^p}^{q}).
\end{align*}
As for $\sI_5(t,\eps)$,  by  Burkholder-Davis-Gundy type inequality and assumption (\ref{asuq}), we obtain
\begin{align*}
\sJ_{5}(t,\eps)&\leq C_5\,\eps^{q}\left(\int_0^t\mE\|(-A)^\gamma e^{(t-s)A}D_x\psi(X_s^\eps,Y_s^{\eps})\Sigma(X_s^\eps,Y_s^{\eps})Q_1^{1/2}\|_{\sR(H,L^p)}^{2}\dif  s\right)^{q/2}\\
&\leq  C_5\,\eps^{q}\left(\int_0^t\|(-A)^\gamma e^{(t-s)A}Q_1^{1/2}\|_{\sR(H,L^p)}^{2}\dif  s\right)^{q/2}\leq  C_5\,\eps^{q},
\end{align*}
and similarly, one can check that
\begin{align*}
\sJ_{6}(t,\eps) \leq  C_6\,\eps^{q/2}.
\end{align*}
Combining the above inequalities, we get
the desired estimate (\ref{stfe}).
\end{proof}

Now, we provide the following uniform estimate for $X_t^\eps$.

\bl\label{la41}

  Let $T>0$, $q\geq1$ and $x,y\in L^p$ with $p\in[2,\infty).$ Then for any $\gamma\in[0,1/2)$, we have
\begin{align}\label{msx'}
	\sup\limits_{\eps\in(0,1)}\mE\| (-A)^\gamma X^{\eps}_t\|_{L^p}^{q}\leq C_{\gamma,p,q,T}\,t^{-\gamma q}(1+\|x\|_{L^p}^{q}+\|y\|_{L^p}^{q}),
	\end{align}
where $C_{\gamma,p,q,T}>0$ is a constant.
\el

\begin{proof}
For $\gamma\in[0,1/2),$ by (\ref{spde22}) we have
\begin{align*}
(-A)^\gamma X^{\eps}_t& =(-A)^\gamma e^{tA}x+\int_0^t(-A)^\gamma e^{(t-s)A}F(X^{\eps}_s, Y^{\eps}_s)\dif s\no\\&\quad+\eps^{-1/2}\int_0^t(-A)^\gamma e^{(t-s)A}B(X^{\eps}_s, Y^{\eps}_s)\dif s\no\\
&\quad+\int_0^t(-A)^\gamma e^{(t-s)A}\Sigma(X_s^\eps,Y_s^\eps)\dif W_s^1
=:\sum\limits_{i=1}^4\sX_i(t,\eps).
\end{align*}
For the first term, it follows from (\ref{pp1}) directly that
$$\mE\|\sX_1(t,\eps)\|_{L^p}^q\leq C_1\,t^{-\gamma}\|x\|_{L^p}^q.$$
To control the second term, by Minkowski's inequality, we have
\begin{align*}
\mE\|\sX_2(t,\eps)\|_{L^p}^q&\leq C_2\Big(\int_0^t\big(\mE\|(-A)^\gamma e^{(t-s)A}F(X^{\eps}_s, Y^{\eps}_s)\|_{L^p}^q\big)^{1/q}\dif s\Big)^q\\
&\leq C_2\Big(\int_0^t(t-s)^{-\gamma}\dif s\Big)^q\leq C_2.
\end{align*}
As for $\sX_3(t,\eps),$ since $B(x,y)$ satisfies the centering condition (\ref{mcen}), by applying Lemma \ref{strf}  we have
$$\mE\|\sX_3(t,\eps)\|_{L^p}^q\leq C_3\,t^{-\gamma q} (1+\|x\|_{L^p}^{q}+\|y\|_{L^p}^{q}).$$
Finally, by  Burkholder-Davis-Gundy type inequality and assumption (\ref{asuq}), we deduce that
\begin{align*}
\mE\|\sX_4(t,\eps)\|^q&\leq C_4\Big(\int_0^t\mE\|(-A)^\gamma e^{(t-s)A}\Sigma(X_s^\eps,Y_s^\eps)Q_1^{1/2}\|_{\sR(H,L^p)}^2\dif s\Big)^{q/2}\leq  C_4.
\end{align*}
Combining the above inequalities, we get the desired estimate (\ref{msx'}).
\end{proof}



\section{Diffusion approximation}

\subsection{Kolmogorov equation}

Note that the process $\bar X_t$ depends on the initial value $x$. Below, we shall write $\bar X_t(x)$  when we want to stress its dependence on the
initial value.
Let $\bar \cL$ be  the infinitesimal generator of the  Markov process $\bar X_t$, i.e.,
\begin{align}
\bar \cL\varphi(x):=\bar \cL(x)\varphi(x):=(\bar \cL_0(x)+\bar \cL_1(x))\varphi(x):=(\bar \cL_0+\bar \cL_1)\varphi(x),\quad\forall\varphi\in C_b^2(H),\label{bL3}
\end{align}
where $\bar \cL_0$ and $\bar \cL_1$ are given by
\begin{align}\label{bL0}
\bar \cL_0\varphi(x):=\<\overline{B\cdot\nabla_x\Psi}(x),D_x\varphi(x)\>
+\frac{1}{2}\,Tr\big[D^2_{x}\varphi(x)\Upsilon(x)\Upsilon^*(x)\big]
\end{align}
and
\begin{align}\label{bL1}
\bar \cL_1\varphi(x):=\<Ax+\bar F(x),D_x\varphi(x)\>
+\frac{1}{2}\,Tr\big[D^2_{x}\varphi(x)\bar\Sigma(x)Q_1\bar\Sigma^*(x)\big].
\end{align}
  Fix $T>0$, consider the following Cauchy problem on $[0,T]\times H$:
\begin{equation} \label{kez}
\left\{ \begin{aligned}
&\partial_t\bar u(t,x)=\bar \cL\,\bar u(t,x),\quad t\in(0,T],\\
& \bar u(0,x)=\varphi(x),\\
\end{aligned} \right.
\end{equation}
where $\varphi: H\to\mR$ is measurable.
We have the following result, which will be used below to prove the weak convergence of $X_t^\eps$ to $\bar X_t$.

\bt\label{lako}
For every $\varphi\in C_b^4(H)$, there exists a   solution to equation (\ref{kez}) which is given by
\begin{align}\label{prou}
\bar u(t,x)=\mE\big[\varphi(\bar X_t(x))\big].
\end{align}
Moreover, we have:
	
	\vspace{1mm}
	\noindent
(i) for any  $t\in(0,T]$, $x\in H$ and $h\in \cD((-A)^{\beta})$ with $\beta\in[0,1)$,
	\begin{align}\label{uz}
	|D_x\bar u(t,x).(-A)^\beta h|\leq C_1\,t^{-\beta}(1+\|x\|_{L^4})\| h\|_{L^4};
	\end{align}
(ii) for any $t\in(0,T]$, $x\in H$, $h_1\in \cD((-A)^{\beta_1})$ and $h_2\in \cD((-A)^{\beta_2})$  with $\beta_1,\beta_2\in[0,1/2)$,
	\begin{align}\label{uz2}
	| D^2_{x}\bar u(t,x).((-A)^{\beta_1}h_1,(-A)^{\beta_2}h_2)|\leq C_2\,t^{-\beta_1-\beta_2}(1+\|x\|_{L^4})\| h_1\|_{L^8}\| h_2\|_{L^8};
	\end{align}
(iii) for any $t\in(0,T]$, $x,h_2,h_3\in H$ and $h_1\in \cD((-A)^{\beta_1})$ with $\beta_1\in[0,1/2)$,
	\begin{align}\label{uz3}
	| D_{x}^3\bar u(t,x).((-A)^{\beta_1}h_1,h_2,h_3)|\leq\! C_3\,t^{-\beta_1}\| h_1\|\| h_2\|\|h_3\|;
\end{align}
(iv) for any $t\in(0,T]$ and $x,h_1,h_2,h_3,h_4 \in H,$
	\begin{align}\label{uz4}
	| D_{x}^4\bar u(t,x).(h_1,h_2,h_3,h_4)|\leq\! C_4\,\| h_1\|\| h_2\|\|h_3\|\|h_4\|;
\end{align}
(v) For any $t\in(0,T],x\in \cD((-A)^{\vartheta_1})$ with $\vartheta_1\in(0,1/2)$ and $h\in \cD((-A)^{\vartheta_2})$ with $\vartheta_2\in(0,1/4)$,
	\begin{align}\label{utz}
	\vert \p_tD_x\bar u(t,x).h|\leq C_5\,(1+\|x\|_{L^4})\big(t^{-1+\vartheta_1+\vartheta_2}\|(-A)^{\vartheta_1}x\|_{L^8} +t^{-1+\vartheta_2}\big)\|(-A)^{\vartheta_2}h\|_{L^8};
	\end{align}
(vi)  for any $t\in(0,T], x\in \cD((-A)^{\vartheta_1})$ with $\vartheta_1\in(0,1/2),$ $h_1\in \cD((-A)^{\vartheta_2})$ with $\vartheta_2\in(0,1/4)$ and $h_2\in \cD((-A)^{\vartheta_3})$ with $\vartheta_3\in(0,1/4)$ satisfying $\vartheta_1+\vartheta_2+\vartheta_3>1/2,$
	\begin{align}\label{utxx}
	\vert \p_tD_x^2\bar u(t,x).(h_1,h_2)|\leq  C_6\,\big(t^{-1+\vartheta_1+\vartheta_2+\vartheta_3}&
\|(-A)^{\vartheta_1}x\|+t^{-1+\vartheta_2+\vartheta_3}(1+\|x\|_{L^4})\big)\no\\
&\times\|(-A)^{\vartheta_2}h_1\|_{L^8}\|(-A)^{\vartheta_3}h_2\|_{L^8},
	\end{align}
where $C_i$, $i=1,\cdots, 6,$ are positive constants.
\et

\begin{proof}
The estimates {\it (i)}-{\it (iii)} have been proven in \cite[Theorem 4.2, Theorem 4.3 and Proposition 4.5]{Br4}, while estimate {\it(iv)} follows by (\ref{prou}).

\vspace{1mm}
\noindent {\it (v)}
To  prove estimate (\ref{utz}), by (\ref{kez}) we have that for any $h\in H,$
\begin{align}\label{ptdz}
\partial_tD_x\bar u(t,x).h=D_x\partial_t\bar u(t,x).h=D_x(\bar\cL\bar u(t,x)).h.
\end{align}
By definition (\ref{bL3}),  we get
\begin{align*}
D_x\bar\cL\bar u(t,x).h&=D_x^2\bar u(t,x).(Ax+\bar F(x)+\overline{B\cdot\nabla_x\Psi}(x),h)\no\\
&\quad+\<Ah+D_x\bar F(x).h+D_x(\overline{B\cdot\nabla_x\Psi}(x)).h,D_x\bar u(t,x)\>\no\\
&\quad+\frac{1}{2}\sum\limits_{n=1}^{\infty} D_x^3\bar u(t,x).(\Upsilon(x)e_{n},\Upsilon(x)e_{n},h)\no\\
&\quad+\sum\limits_{n=1}^{\infty} D_x^2\bar u(t,x).((D_x\Upsilon(x).h)e_{n},\Upsilon(x)e_{n})\no\\
&\quad+\frac{1}{2}\sum\limits_{n=1}^{\infty} D_x^3\bar u(t,x).(\bar\Sigma(x)Q_1^{1/2}e_{n},\bar\Sigma(x)Q_1^{1/2}e_{n},h)\no\\
&\quad+\sum\limits_{n=1}^{\infty} D_x^2\bar u(t,x).((D_x\bar\Sigma(x).h)Q_1^{1/2}e_{n},\bar\Sigma(x)Q_1^{1/2}e_{n}).
\end{align*}
Recall that
$$D_x^2\bar u(t,x).(v_1,v_2)=D_x^2\bar u(t,x).(v_2,v_1),\forall v_1,v_2\in H,$$
i.e., $D_x^2\bar u(t,x)\in\sL(H)$ is self-adjoint. As a result, we have that for $\gamma\in[0,1],$
\begin{align}\label{lzdz'}
D_x^2\bar u(t,x).(Av_1,v_2)&=\<D_x^2\bar u(t,x).v_2,Av_1\>\no\\
&=\<(-A)^\gamma D_x^2\bar u(t,x).v_2,(-A)^{1-\gamma}v_1\>\no\\
&=\< D_x^2\bar u(t,x).(-A)^\gamma v_2,(-A)^{1-\gamma}v_1\>\no\\
&=D_x^2\bar u(t,x).((-A)^{1-\gamma}v_1,(-A)^{\gamma}v_2).
\end{align}
Using estimates (\ref{uz}), (\ref{uz2}), (\ref{uz3}) and (\ref{lzdz'}), we deduce that for $\vartheta_1\in(0,1/2)$ and $\vartheta_2\in(0,1/4),$
\begin{align}\label{lzdz}
&|D_x\bar\cL\bar u(t,x).h|\no\\
&\leq C_1\|h\|_{\infty}+D_x^2\bar u(t,x).(Ax,h)+\<Ah,D_x\bar u(t,x)\>\no\\
&= C_1\|h\|_{\infty}+D_x^2\bar u(t,x).((-A)^{1/2-\vartheta_1}(-A)^{\vartheta_1}x,(-A)^{1/2-\vartheta_2}(-A)^{\vartheta_2} h)\no\\&\quad+\<(-A)^{1-\vartheta_2}(-A)^{\vartheta_2}h,D_x\bar u(t,x)\>\no\\
&\leq  C_1\,(1+\|x\|_{L^4})\big(t^{-1+\vartheta_1+\vartheta_2}\|(-A)^{\vartheta_1}x\|_{L^8} +t^{-1+\vartheta_2}\big)\|(-A)^{\vartheta_2}h\|_{L^8},
\end{align}
where in the last inequality, we also used the Sobolev inequality that $\|h\|_{\infty}\leq c_0\|(-A)^{\vartheta_2}h\|_{L^8}.$
Combining (\ref{ptdz}) and (\ref{lzdz}), we obtain (\ref{utz}).

\vspace{1mm}
\noindent {\it (vi)} In view of (\ref{ptdz}), we note that for any $h_1,h_2\in H,$
\begin{align*}
\partial_tD_x^2\bar u(t,x).(h_1,h_2)&=D_x^3\bar u(t,x).(Ax+\bar F(x)+\overline{B\cdot\nabla_x\Psi}(x),h_1,h_2)\\
&\quad+D_x^2\bar u(t,x).(Ah_2+D_x\bar F(x).h_2+D_x(\overline{B\cdot\nabla_x\Psi}(x)).h_2,h_1)\\
&\quad+D_x^2\bar u(t,x).(Ah_1+D_x\bar F(x).h_1+D_x(\overline{B\cdot\nabla_x\Psi}(x)).h_1,h_2)\\
&\quad+\<D_x^2\bar F(x).(h_1,h_2)+D_x^2(\overline{B\cdot\nabla_x\Psi}(x)).(h_1,h_2),D_x\bar u(t,x)\>\\
&\quad+\frac{1}{2}\sum\limits_{n=1}^{\infty} D_x^4\bar u(t,x).(\Upsilon(x)e_{n},\Upsilon(x)e_{n},h_1,h_2)\no\\
&\quad+\sum\limits_{n=1}^{\infty} D_x^3\bar u(t,x).((D_x\Upsilon(x).h_2)e_{n},\Upsilon(x)e_{n},h_1)\no\\
&\quad+\sum\limits_{n=1}^{\infty} D_x^3\bar u(t,x).((D_x\Upsilon(x).h_1)e_{n},\Upsilon(x)e_{n},h_2)\no\\
&\quad+2\sum\limits_{n=1}^{\infty} D_x^2\bar u(t,x).((D_x^2\Upsilon(x).(h_1,h_2))e_{n},\Upsilon(x)e_{n})\no\\
&\quad+\frac{1}{2}\sum\limits_{n=1}^{\infty} D_x^4\bar u(t,x).(\bar\Sigma(x)Q_1^{1/2}e_{n},\bar\Sigma(x)Q_1^{1/2}e_{n},h_1,h_2)\no\\
&\quad+\sum\limits_{n=1}^{\infty} D_x^3\bar u(t,x).((D_x\bar\Sigma(x).h_2)Q_1^{1/2}e_{n},\bar\Sigma(x)Q_1^{1/2}e_{n},h_1)\\
&\quad+\sum\limits_{n=1}^{\infty} D_x^3\bar u(t,x).((D_x\bar\Sigma(x).h_1)Q_1^{1/2}e_{n},\bar\Sigma(x)Q_1^{1/2}e_{n},h_2)\\
&\quad+2\sum\limits_{n=1}^{\infty} D_x^2\bar u(t,x).((D_x^2\bar\Sigma(x).(h_1,h_2))Q_1^{1/2}e_{n},\bar\Sigma(x)Q_1^{1/2}e_{n}).
\end{align*}
For $\gamma_1,\gamma_2\in[0,1],$ by applying (\ref{lzdz'}) we have
\begin{align}\label{lzdz''}
D_x^3\bar u(t,x).(Av_1,v_2,v_3)&=D_x\<D_x^2\bar u(t,x).v_2,Av_1\>.v_3\no\\
&=D_x\< D_x^2\bar u(t,x).(-A)^\gamma_1v_2,(-A)^{1-\gamma_1}v_1\>.v_3\no\\
&=D_x^3\bar u(t,x).((-A)^{1-\gamma_1}v_1,(-A)^{\gamma_1}v_2,v_3)\no\\
&=D_x^3\bar u(t,x).((-A)^{1-\gamma_1}v_1,v_3,(-A)^{\gamma_1}v_2)\no\\
&=D_x^3\bar u(t,x).((-A)^{1-\gamma_1-\gamma_2}v_1,(-A)^{\gamma_2}v_3,(-A)^{\gamma_1}v_2)\no\\
&=D_x^3\bar u(t,x).((-A)^{1-\gamma_1-\gamma_2}v_1,(-A)^{\gamma_1}v_2,,(-A)^{\gamma_2}v_3).
\end{align}
Using (\ref{uz})-(\ref{uz4}) and (\ref{lzdz''}), one can check that
\begin{align*}
&|\partial_tD_x^2\bar u(t,x).(h_1,h_2)|\\
&
\leq C_2\|h_1\|_{\infty}\|h_2\|_{\infty}+D_x^3\bar u(t,x).(Ax,h_1,h_2)+D_x^2\bar u(t,x).(Ah_1,h_2)\\&\quad+D_x^2\bar u(t,x).(Ah_2,h_1)\\
&=C_2\|h_1\|_{\infty}\|h_2\|_{L^\infty}+D_x^3\bar u(t,x).((-A)^{1-\vartheta_1-\vartheta_2-\vartheta_3}(-A)^{\vartheta_1}x,(-A)^{\vartheta_2}h_1,
(-A)^{\vartheta_3}h_2)\\
&\quad+2D_x^2\bar u(t,x).((-A)^{1/2-\vartheta_2}(-A)^{\vartheta_2}h_1,(-A)^{1/2-\vartheta_3}
(-A)^{\vartheta_3}h_2)\\
&\leq C_2\,(t^{-1+\vartheta_1+\vartheta_2+\vartheta_3}
\|x\|_{(-A)^{\vartheta_1}}+t^{-1+\vartheta_2+\vartheta_3}(1+\|x\|_{L^4}))
\|(-A)^{\vartheta_2}h_1\|_{L^8}\|(-A)^{\vartheta_3}h_2\|_{L^8}.
\end{align*}
The proof is finished.
\end{proof}

\subsection{Proof of Theorem \ref{main3}}
The following weak fluctuation estimates for an integral functional of $(X_t^\eps,Y_t^{\eps})$  will play an important role in proving (\ref{nzz}).

\bl[Weak fluctuation estimates]\label{weaf}
Let $T>0$ and $x,y\in L^8$.  Then,

\vspace{1mm}
\noindent (i)
for any $\phi(\cdot,\cdot): H\times H \to H$ satisfying both (i)-(ii) of Lemma \ref{fbg} and the centering condition (\ref{cen}),
 we have
\begin{align}
\bigg|\mE\left(\int_0^T\<\phi(X_t^{\eps},Y_t^{\eps}),D_x\bar u(T-t,X_t^\eps)\>\dif t\right)\bigg|&\leq C_T\,\eps^{\frac{1}{2}};  \label{we1}
\end{align}

\vspace{1mm}
\noindent (ii)
 for any $\tilde \phi(\cdot,\cdot):H\times H \to\sL(H)$ satisfying both Lemma \ref{sig} and the  centering condition
\begin{align}\label{cen222}
\int_{H}Tr[D_x^2\bar u(T-t,x)\tilde \phi(x,y)]\mu^x(\dif y)=0, \end{align}
 we have
\begin{align}
\bigg|\mE\left(\int_0^T Tr[D_x^2\bar u(T-t,X_t^\eps)\tilde \phi(X_t^\eps,Y_t^\eps)]\dif t\right)\bigg|&\leq C_T\,\eps^{\frac{1}{2}};  \label{we2}
\end{align}
where $C_T>0$ is a constant.
\el

\begin{proof}

\vspace{1mm}
\noindent {\it (i)} Let $\psi(x,y)$ solve the Poisson equation $$
\cL_2(x,y)\psi(x,y)=-\phi(x,y),
$$
 and define \begin{align}\label{defp}\psi_t(x,y)=\<\psi(x,y),D_x\bar u(T-t,x)\>.\end{align}
 It is easy to check that $\psi_t(x,y)$ solves the following Poisson equation:
\begin{align}\label{psi1}
\cL_2(x,y)\psi_t(x,y)=-\<\phi(x,y),D_x\bar u(T-t,x)\>.
\end{align}
  According to Theorems \ref{PP} and \ref{lako}, we can apply It\^o's formula to $\psi_t(X_t^\eps,Y_t^\eps)$ to derive that
	\begin{align}\label{aa}
	\mE[\psi_T(X_T^{\eps},Y_T^\eps)]&=\psi_0(x,y) +\eps^{-1}\mE\left(\int_0^T\mathcal{L}_2\psi_t(X_t^{\eps},Y_t^\eps)\dif t\right)\no\\
	&\quad+\mE\left(\int_0^T(\partial_t+\eps^{-1/2}\cL_0+\cL_1)\psi_t(X_t^{\eps},Y_t^\eps)\dif t\right),
	\end{align}
where $\cL_0,\cL_1$ and $\cL_2$ are defined by (\ref{L0}), (\ref{L1}) and (\ref{L2}), respectively.
	Multiplying  both sides of (\ref{aa}) by $\eps$ and taking into account (\ref{psi1}), we get
	\begin{align*}
	&\bigg|\mE\left(\int_0^T\<\phi(X_t^{\eps},Y_t^{\eps}),D_x\tilde u(t,X_t^\eps)\>\dif t\right)\bigg|\\&=\bigg|\mE\left(\int_0^T\mathcal{L}_2\psi_t(X_t^{\eps},Y_t^\eps)\dif t\right)\bigg|\\ &\leq\eps\,\big|\mE\big[\psi_0(x,y)-\psi_T(X_T^{\eps},Y_T^\eps)\big]\big|+\eps\,\bigg|\mE\left(\int_0^T
\partial_t\psi_t(X_t^{\eps},Y_t^{\eps})\dif t\right)\bigg|\\
&\quad+\sqrt{\eps}\,\bigg|\mE\left(\int_0^T\cL_0\psi_t(X_t^{\eps},Y_t^{\eps})\dif t\right)\bigg|+\eps\,\bigg|\mE\left(\int_0^T\cL_1\psi_t(X_t^{\eps},Y_t^{\eps})\dif t\right)\bigg|=:\sum_{i=1}^4\sR_i(T,\eps).
	\end{align*}
By applying (\ref{defp}), (\ref{uz}) and Theorem \ref{PP}, we have
\begin{align*}
\sR_1(T,\eps)\leq C_1\,\eps\big(\|\psi(x,y)\|+\|\psi(X_T^{\eps},Y_T^{\eps})\|\big)
\leq C_1\,\eps.
\end{align*}
 To control the second term, by (\ref{defp}), (\ref{utz}), (\ref{lfp}), (\ref{msy}) and Lemma \ref{la41}, we have for any $\vartheta_1\in(0,1/2),$ $\vartheta_2\in(0,1/4)$ and small enough $\delta>0,$
\begin{align*}
&\sR_2(T,\eps)\leq\eps\,\mE\bigg|\int_0^T\<\psi(X_t^{\eps},Y_t^{\eps}),\partial_tD_x\bar u(T-t,X_t^\eps)\>\dif t\bigg|\\
&\leq C_2\,\eps\mE\bigg(\int_0^T (1+\|X_t^\eps\|_{L^4})\\&\quad\times\big((T-t)^{-1+\vartheta_1+\vartheta_2}\|(-A)^{\vartheta_1} X_t^\eps\|_{L^8}+(T-t)^{-1+\vartheta_2}\big)\|(-A)^{\vartheta_2} \psi(X_t^\eps,Y_t^\eps)\|_{L^8}\dif t\bigg)\\
&\leq C_2\,\eps\int_0^T (T-t)^{-1+\vartheta_2}\big(\mE(1+\|X_t^\eps\|_{L^4})^2\big)^{1/2}\\
&\qquad\quad\qquad\times\big(\mE(1+
\|(-A)^{\vartheta_1} X_t^\eps\|_{L^8}^2+\|(-A)^{\vartheta_2+\delta} Y_t^\eps\|_{L^8}^2)^2\big)^{1/2}\dif t\\
&\leq  C_2\,\eps\int_0^T(T-t)^{-1+\vartheta_2}(t^{-2\vartheta_1}+t^{-2\vartheta_2-2\delta})\dif t\leq C_2\,\eps.
\end{align*}
For the third term, by definitions (\ref{L0}), (\ref{defp}), Theorems \ref{PP} and \ref{lako},  we have
\begin{align*}
\sR_3(T,\eps)&\leq\sqrt{\eps}\,\mE\bigg|\int_0^T\< B(X_t^\eps,Y_t^\eps),D_x\psi_t(X_t^\eps,Y_t^\eps)\>\dif t\bigg|\\
&\leq\sqrt{\eps}\,\mE\bigg|\int_0^T\< D_x\psi(X_t^\eps,Y_t^\eps).B(X_t^\eps,Y_t^\eps),D_x\bar u(T-t,X_t^\eps)\>\dif t\bigg|\\
&\quad+\sqrt{\eps}\,\mE\bigg|\int_0^TD_x^2\bar u(T-t,X_s^\eps).(B(X_t^\eps,Y_t^\eps),\psi(X_t^\eps,Y_t^\eps))\dif t\bigg|\\
&\leq C_3\,\eps^{1/2}.
\end{align*}
For the last term, it is easy to check that
\begin{align*}
\sR_{4}(T,\eps)&\leq\eps\,\mE\bigg|\int_0^T\<AX_t^\eps,D_x\psi_t(X_t^\eps,Y_t^\eps)\>\dif t\bigg|\\&\quad+\eps\,\mE\bigg|\int_0^TF(X_t^\eps,Y_t^\eps),D_x\psi_t(X_t^\eps,Y_t^\eps)\>\dif t\bigg|\\
&\quad+\frac{\eps}{2}\,\mE\bigg|\int_0^TTr\big[D_x^2\psi_t(X_t^\eps,Y_t^\eps)
\Sigma(X_t^\eps,Y_t^\eps)Q_1\Sigma^*(X_t^\eps,Y_t^\eps)\big]
\dif t\bigg|\\&\leq \eps\,\mE\bigg|\int_0^T\<AX_t^\eps,D_x\psi_t(X_t^\eps,Y_t^\eps)\>\dif t\bigg|+C_4\,\eps=:\sR_{4,1}(T,\eps)+C_4\,\eps.
\end{align*}
In view of (\ref{defp}), (\ref{lzdz'}), (\ref{uz}), (\ref{uz2}) and (\ref{lfp}),  we have for any $x\in \cD((-A)^{\vartheta_1})$  with $\vartheta_1\in(0,1/2)$, $y\in \cD((-A)^{\vartheta_2})$ with $\vartheta_2\in(0,1/4)$ and small enough $\delta>0,$
\begin{align*}
	|\<Ax,D_x\psi_t(x,y)\>|&= D^2_{x}\bar u(T-t,x).(\psi(x,y),A x)+\<D_x\psi(x,y).A x,D_x\bar u(T-t,x)\>\\
&= D^2_{x}\bar u(T-t,x).((-A)^{1/2-\vartheta_2}(-A)^{\vartheta_2} \psi(x,y),(-A)^{1/2-\vartheta_1}(-A)^{\vartheta_1}  x)\\&\quad+\< D_x\psi(x,y). (-A)^{\vartheta_1} x, (-A)^{1-\vartheta_1} D_x\bar u(T-t,x)\>\\
&\leq C_4\,(T-t)^{-1+\vartheta_1+\vartheta_2}(1+\|x\|_{L^4})\|(-A)^{\vartheta_2} \psi(x,y)\|_{L^8}\|(-A)^{\vartheta_1}  x\|_{L^8}\\
&\quad+C_4\,(T-t)^{-1+\vartheta_1}(1+\|x\|_{L^4})\|D_x\psi(x,y). (-A)^{\vartheta_1} x\|_{L^4}\\
&\leq  C_4\,(T-t)^{-1+\vartheta_1}(1+\|x\|_{L^4})\big(1+\|(-A)^{\vartheta_1}x\|_{L^8}^2
+\|(-A)^{\vartheta_2+\delta}y\|_{L^8}^2\big).
	\end{align*}
Consequently, by  Lemmas {\ref{la42}} and \ref{la41}  we have
\begin{align*}
\sR_{4,1}(T,\eps)&\leq C_4\,\eps\int_0^T (T-t)^{-1+\vartheta_1}\big(\mE(1+\|X_t^\eps\|_{L^4})^2\big)^{1/2}\\
&\qquad\quad\qquad\qquad\qquad\times\big(\mE(1+
\|(-A)^{\vartheta_1} X_t^\eps\|_{L^8}^2+\|(-A)^{\vartheta_2+\delta} Y_t^\eps\|_{L^8}^2)^2\big)^{1/2}\dif t\\
&\leq C_4\,\eps \int_0^T(T-t)^{-1+\vartheta_1}(t^{-2\vartheta_1}+t^{-2\vartheta_2-2\delta})\dif t\leq C_4\, \eps.
\end{align*}
Combining the above inequalities, we get estimate (\ref{we1}).

\vspace{1mm}
\noindent {\it (ii)}
Consider the following Poisson equation:
\begin{align}\label{psi2}
\cL_2(x,y)\tilde\psi_t(x,y)=-Tr[D_x^2\bar u(T-t,x)\tilde \phi(x,y)]=:-\tilde\phi_t(x,y).
\end{align}
Using exactly the same arguments as above, we can obtain
	\begin{align*}
	&\bigg|\mE\left(\int_0^T\tilde \phi_t(X_t^\eps,Y_t^\eps)\dif t\right)\bigg|=\bigg|\mE\left(\int_0^T\mathcal{L}_2\tilde\psi_t(X_t^{\eps},Y_t^\eps)\dif t\right)\bigg|\\ &\leq\eps\,\mE\big|\big[\tilde\psi_0(x,y)-\tilde\psi_T(X_T^{\eps},Y_T^\eps)\big]\big|
+\eps\,\mE\bigg|\int_0^T\partial_t\tilde\psi_t(X_t^{\eps},Y_t^{\eps})\dif t\bigg|\\
&\quad+\sqrt{\eps}\,\mE\bigg|\int_0^T\cL_0\tilde\psi_t(X_t^{\eps},Y_t^{\eps})\dif t\bigg|+\eps\,\mE\bigg|\int_0^T\cL_1\tilde\psi_t(X_t^{\eps},Y_t^{\eps})\dif t\bigg|=:\sum_{i=1}^4\sV_i(T,\eps).
	\end{align*}
According to definitions (\ref{L0}), (\ref{psi2}), Theorems \ref{PP} and \ref{lako}, it is easy to check that
\begin{align*}
\sV_1(T,\eps)+\sV_3(T,\eps)
\leq C_1\,\eps^{1/2}.
\end{align*}
 To estimate the second term, by making use of (\ref{utxx}), (\ref{uz3}) and (\ref{lfp}),
we have that for any $x\in \cD((-A)^{\vartheta_1})$  with $\vartheta_1\in(0,1/2)$, $y\in \cD((-A)^{\vartheta_2})$ with $\vartheta_2\in(0,1/4)$ satisfying $\vartheta_1+2\vartheta_2>1/2$ and small enough $\delta>0,$
\begin{align*}
	|\p_t\tilde \phi_t(x,y)|&\leq\Big|\sum\limits_{n=0}^{\infty}\p_tD_x^2\bar u(T-t,x).((\tilde\phi^{1/2}(x,y)
e_n,(\tilde\phi^{1/2}(x,y))^*e_n)\Big|\\
&\leq C_2\,\big((T-t)^{-1+2\vartheta_2}(1+\|x\|_{L^4})+(T-t)^{-1+\vartheta_1+2\vartheta_2}
\|(-A)^{\vartheta_1}x\|\big)\\&\quad\times(1+\|(-A)^{\vartheta_2+\delta}x\|_{L^8}^2
+\|(-A)^{\vartheta_2
+\delta}y\|_{L^8}^2).
	\end{align*}
Thus, by definition (\ref{psi2}), Theorem \ref{PP}, Lemmas \ref{la42} and \ref{la41}, we deduce that
\begin{align*}
\sV_2(T,\eps)
\leq C_2\,\eps\int_0^T(T-t)^{-1+2\vartheta_2}t^{-\vartheta_1-2\vartheta_2-2\delta}\dif t\leq C_2\,\eps.
\end{align*}
For the last term, we have
\begin{align*}
\sV_{4}(T,\eps)&\leq C_0\,\eps+\eps\,\mE\bigg|\int_0^T\<AX_t^\eps,D_x\tilde\psi_t(X_t^\eps,Y_t^\eps)\>\dif t\bigg|=:C_0\,\eps+\sV_{4,1}(T,\eps).
\end{align*}
As for $\sV_{4,1}(T,\eps)$, by (\ref{lzdz'}), (\ref{lzdz''}), (\ref{uz2}) and (\ref{uz3}) we have for any $x\in \cD((-A)^{\vartheta_1})$  with $\vartheta_1\in(0,1/2)$ and $y\in \cD((-A)^{\vartheta_2})$ with $\vartheta_2\in(0,1/4)$ satisfying $\vartheta_1+\vartheta_2>1/2,$ and for small enough $\delta>0,$
	\begin{align*}
	&|\<Ax,D_x\tilde\phi_t(x,y)\>|\\&\leq \Big|\sum\limits_{n=0}^{\infty}D^3_{x}\bar u(T-t,x).((-A)^{\vartheta_2}\tilde\phi^{1/2}(x,y)e_n,(\tilde\phi^{1/2}(x,y))^*e_n,(-A)^{1-
\vartheta_1-\vartheta_2}(-A)^{\vartheta_1} x)\Big|\\&+2\Big|\sum\limits_{n=0}^{\infty}D^2_{x}\bar u(T-t,x).(D_x\tilde\phi^{1/2}(x,y).((-A)^{\vartheta_1} x)
e_n,(-A)^{1-\vartheta_1-\vartheta_2}(-A)^{\vartheta_2}(\tilde\phi^{1/2}(x,y))^*e_n)\Big|\\
&\leq  C_4\,(T-t)^{-1+\vartheta_1+\vartheta_2}(1+\|x\|_{L^4})(1+\|(-A)^{\vartheta_1}x\|_{L^8}^2
+\|(-A)^{\vartheta_2+\delta}y\|_{L^8}^2).
	\end{align*}
Thus by Theorem \ref{PP}, Lemmas \ref{la42} and \ref{la41}, we have
\begin{align*}
\sV_{4}(T,\eps)
&\leq C_0\,\eps+C_4\,\eps\int_0^T(T-t)^{-1+\vartheta_1+\vartheta_2}(t^{-2\vartheta_1}
+t^{-2\vartheta_2-2\delta})\dif t\leq C_4\,\eps.
	\end{align*}
Combining the above inequalities, we get estimate (\ref{we2}).
\end{proof}

Now, we are in the position to give:

\begin{proof}[Proof of Theorem \ref{main3}]
	Given $T>0$ and $\varphi\in C_b^4(H)$, let $\bar u$ solve the Cauchy problem (\ref{kez}). For $t\in[0,T]$ and $x\in H$, define
$$
\tilde u(t,x)=\bar u(T-t,x).
$$
Then one can check that
$$
\tilde u(T,x)=\bar u(0,x)=\varphi(x)\quad\text{and}\quad\tilde u(0,x)=\bar u(T,x)=\mE[\varphi(\bar X_T(x))].
$$
Using It\^o's formula and by definitions (\ref{L}), (\ref{bL3}), (\ref{bL0}), (\ref{bL1}) and equality (\ref{kez}), we deduce that
\begin{align*}
\mE[\varphi(X_T^{\eps})]-\mE[\varphi&(\bar X_T)]=\mE[\tilde u(T,X_T^{\eps})-\tilde u(0,x)]=\mE\left(\int_0^T\big(\p_t+\cL\big)\tilde u(t,X_t^{\eps})\dif t\right)\no\\
&=\mE\left(\int_0^T[\mathcal{L}\tilde u(t,X_t^{\eps})-\mathcal{\bar L}\tilde u(t,X_t^{\eps})]\dif t\right)\no\\
&=\mE\left(\int_0^T(\mathcal{L}_1-\mathcal{\bar L}_1)\tilde u(t,X_t^{\eps})\dif t\right)+\mE\left(\int_0^T(\eps^{-1/2}\mathcal{L}_0-\mathcal{\bar L}_0)\tilde u(t,X_t^{\eps})\dif t\right)\no\\
&=\mE\left(\int_0^T\langle  F(X_t^\eps,Y_t^\eps)-\bar F(X_t^\eps), D_x\tilde u(t,X_t^{\eps})\rangle \dif t\right)\no\\
&\quad+\frac{1}{2}\mE\bigg(\int_0^T\big(Tr\big[D_x^2\tilde u(t,X_t^{\eps})\Sigma(X_t^\eps,Y_t^\eps)Q_1\Sigma^*(X_t^\eps,Y_t^\eps)\big]\no\\&
\qquad\qquad\qquad\qquad -Tr\big[D_x^2\tilde u(t,X_t^{\eps})\bar\Sigma(X_t^\eps)Q_1\bar\Sigma^*(X_t^\eps)\big]\big) \dif t\bigg)\no\\
&\quad+\mE\bigg(\int_0^T\big(\langle  \eps^{-1/2}B(X_t^\eps,Y_t^\eps)-\overline{B\cdot\nabla_x\Psi}(X_t^\eps), D_x\tilde u(t,X_t^{\eps})\rangle\no\\&
\qquad\qquad\qquad\qquad-\frac{1}{2}\,Tr\big[D_x^2\tilde u(t,X_t^{\eps})\Upsilon(X_t^\eps)\Upsilon^*(X_t^\eps)\big] \dif t\bigg)\no\\
&=:\sum\limits_{i=1}^3\sN_i(T,\eps).
\end{align*}
For the first term, define $$\phi(x,y):=F(x,y)-\bar F(x).$$
It is obvious that $\phi$ satisfies the centering condition (\ref{cen}). 	
Thus, according to (\ref{we1}), we have
$$
|\sN_1(T,\eps)|\leq C_1\,\eps^{1/2}.
$$
To control the second term, let
\begin{align*}
\sI_t(x,y)=:Tr[D_x^2\tilde u(t,x)(\Sigma(x,y)Q_1\Sigma^*(x,y)
-\bar\Sigma(x)Q_1\bar\Sigma^*(x))].
\end{align*}
In view of (\ref{df2}), one can check that $\sI_t(x,y)$ satisfies the centering condition (\ref{cen222}).
Thus, by applying (\ref{we2})  we have
$$|\sN_2(T,\eps)|\leq C_2\,\eps^{1/2}.$$
For the last term, recall that $\Psi$  solves the Poisson equation (\ref{poF}), and define
	$$
	\tilde\Psi_t(x,y):=\<\Psi(x,y),D_x\tilde u(t,x)\>.
	$$
	Since $\cL_2$ is an operator with respect to the $y$-variable, one can check that $\tilde\Psi_t$ solves the following Poisson equation:
	\begin{align*}
	\mathcal{L}_2(x,y)\tilde\Psi_t(x,y)=-\langle B(x,y),D_x\tilde u(t,x).
	\end{align*}
By using exactly the same arguments as in Lemma \ref{weaf}, we can obtain
\begin{align*}
&|\sN_3(T,\eps)|\\&\leq\sqrt{\eps}\,\big|\mE\big[\tilde\Psi_0(x,y)
-\tilde\Psi_T(X_T^{\eps},Y_T^\eps)\big]\big|
+\sqrt{\eps}\,\bigg|\mE\left(\int_0^T\partial_t\tilde\Psi_t(X_t^{\eps},Y_t^{\eps})\dif t\right)\bigg|\\
&\quad+\sqrt{\eps}\,\bigg|\mE\left(\int_0^T\cL_1\tilde\Psi_t(X_t^{\eps},Y_t^{\eps})\dif t\right)\bigg|+\bigg|\mE\left(\int_0^T\cL_0\tilde\Psi_t(X_t^{\eps},Y_t^{\eps})\dif t\right)\\&\quad-\mE\bigg(\int_0^T\big(\<\overline{B\cdot\nabla_x\Psi}(X_t^\eps), D_x\tilde u(t,X_t^{\eps})\>\dif t\bigg)\\&\quad-\mE\bigg(\int_0^T \frac{1}{2}\,Tr\big[D_x^2\tilde u(t,X_t^{\eps})\Upsilon(X_t^\eps)\Upsilon^*(X_t^\eps)\big] \dif t\bigg)\bigg|\\
&\leq C_3\;\eps^{1/2}+\bigg|\mE\left(\int_0^T\< D_x\Psi(X_t^\eps,Y_t^\eps).B(X_t^\eps,Y_t^\eps)-\overline{B\cdot\nabla_x\Psi}(X_t^\eps),D_x\tilde u(t,X_t^\eps)\>\dif t\right)\bigg|\\&
\quad+\bigg|\mE\left(\int_0^TD_x^2\tilde u(t,X_t^\eps).(B(X_t^\eps,Y_t^\eps),\Psi(X_t^\eps,Y_t^\eps))-\frac{1}{2}\,Tr\big[D_x^2\tilde u(t,X_t^{\eps})\Upsilon(X_t^\eps)\Upsilon^*(X_t^\eps)\big]\dif t\right)\bigg|
\\&=:C_3\,\eps^{1/2}+\sN_{3,1}(T,\eps)+\sN_{3,2}(T,\eps).
\end{align*}
For the term $\sN_{3,1}(T,\eps),$ let
\begin{align*}
\sX(x,y)=: D_x\Psi(x,y).B(x,y)-\overline{B\cdot\nabla_x\Psi}(x).
\end{align*}
Note that by the definition of $\overline{B\cdot\nabla_x\Psi}$, one can check that $\sX(x,y)$ satisfies the centering condition (\ref{cen}). Thus, using (\ref{we1}) directly  we obtain
$$\sN_{3,1}(T,\eps)\leq C_3\;\eps^{1/2}.$$
To control the term $\sN_{3,2}(T,\eps),$ let
\begin{align*}
\sY_t(x,y)&=:D_x^2\tilde u(t,x).(B(x,y),\Psi(x,y))-\frac{1}{2}\,Tr\big[D_x^2\tilde u(t,x)\Upsilon(x)\Upsilon^*(x)\big].
\end{align*}
By the definition of $\Upsilon$ in (\ref{sst}), we find that $\sY_t(x,y)$ satisfies the centering condition. As a result of (\ref{we2})  we have
	$$
	\sN_{3,2}(T,\eps)\leq C_3\; \eps^{1/2}.
	$$
Combining the above estimates, we get the desired result.
\end{proof}

\bigskip
\noindent{\bf{Acknowledgements}} The authors are very grateful to the referees for their quite useful suggestions.

\end{document}